%
%
%
\magnification=\magstep1
\newif\ifsect
\secttrue
\def\ssect #1. {\bigbreak\indent{\bf #1.}\enspace\message{#1}}
\def\smallsect #1. #2\par{\bigbreak\noindent{\bf #1.}\enspace{\bf #2}\par
        \global\parano=#1\global\eqnumbo=1\global\thmno=1
        \nobreak\smallskip\nobreak\noindent\message{#2}}
\def\thm #1: #2{\medbreak\noindent{\bf #1:}\if(#2\thmp\else\thmn#2\fi}
\def\thmp #1) { (#1)\thmn{}}
\def\thmn#1#2\par{\enspace{\sl #1#2}\par
        \ifdim\lastskip<\medskipamount \removelastskip\penalty 55\medskip\fi}

\def\qedn{\thinspace\null\nobreak\hfill\hbox{\vbox{\kern-.2pt\hrule height.2pt
depth.2pt\kern-.2pt\kern-.2pt \hbox to2.5mm{\kern-.2pt\vrule width.4pt
\kern-.2pt\raise2.5mm\vbox to.2pt{}\lower0pt\vtop to.2pt{}\hfil\kern-.2pt
\vrule width.4pt\kern-.2pt}\kern-.2pt\kern-.2pt\hrule height.2pt depth.2pt
\kern-.2pt}}\par\medbreak}
\def\pf{\ifdim\lastskip<\smallskipamount \removelastskip\smallskip\fi
        \noindent{\sl Proof\/}:\enspace}
\def\itm#1{\par\indent\llap{\rm #1\enspace}\ignorespaces}
\let\smbar=\bar
\def\bar#1{\overline{#1}}
\def\forclose#1{\hfil\llap{$#1$}\hfilneg}
\def\newforclose#1{
        \ifsect\xdef #1{(\number\parano.\number\eqnumbo)}\else
        \xdef #1{(\number\eqnumbo)}\fi
        \hfil\llap{$#1$}\hfilneg
        \global \advance \eqnumbo by 1}
\def\forevery#1#2$${\displaylines{\let\eqno=\forclose
        \let\neweq=\newforclose\hfilneg\rlap{$\qquad\forall#1$}\hfil#2\cr}$$}
\newcount\parano
\newcount\eqnumbo
\newcount\thmno
\def\neweqt#1$${\xdef #1{(\number\parano.\number\eqnumbo)}
        \eqno #1$$
        \global \advance \eqnumbo by 1}
\def\newthmt#1 #2: #3{\xdef #2{\number\parano.\number\thmno}
        \global \advance \thmno by 1
        \medbreak\noindent{\bf #1 #2:}\if(#3\thmp\else\thmn#3\fi}
\def\neweqf#1$${\xdef #1{(\number\eqnumbo)}
        \eqno #1$$
        \global \advance \eqnumbo by 1}
\def\newthmf#1 #2: #3{\xdef #2{\number\thmno}
        \global \advance \thmno by 1
        \medbreak\noindent{\bf #1 #2:}\if(#3\thmp\else\thmn#3\fi}
\def\bititolo{\empty}

\font\titlefont=cmssbx10 scaled \magstep1
\font\bigfont=cmr10 scaled \magstep1

\font\eightrm=cmr8

\font\ss=cmss10

\nopagenumbers
\binoppenalty=10000
\relpenalty=10000
\let\de=\partial
\def\eps{\varepsilon}
\def\phe{\varphi}

\def\Re{\mathop{\rm Re}\nolimits}

\def\cancel#1#2{\ooalign{$\hfil#1/\hfil$\crcr$#1#2$}}
\def\void{\mathord{\mathpalette\cancel{\mathrel{\hbox{\kern0pt\raise0.8pt\hbox
        {$\scriptstyle\bigcirc$}}}}}}
\def\R{{\rm I\kern-.185em R}}
\def\P{{\rm I\kern-.185em P}}
\def\N{{\rm I\kern-.185em N}}
\def\Z{{\mathchoice{\hbox{\ss Z\kern-.38em Z}}{\hbox{\ss Z\kern-.37em Z}}{\hbox
{\sevenrm Z\kern-.39em Z}}{\hbox{\sevenrm Z\kern-.39em Z}}}}
\def\Q{{\mathchoice{\hbox{\rm\kern.37em\vrule height1.4ex width.05em
depth-.011em\kern-.37em Q}}{\hbox{\rm\kern.37em\vrule height1.4ex width.05em
depth-.011em\kern-.37em Q}}{\hbox{\sevenrm\kern.37em\vrule height1.3ex
width.05em depth-.02em\kern-.3em Q}}{\hbox{\sevenrm\kern.37em\vrule height1.3ex
width.05em depth-.02em\kern-.3em Q}}}}
\def\C{{\mathchoice{\hbox{\rm\kern.37em\vrule height1.4ex width.05em
depth-.011em\kern-.37em C}}{\hbox{\rm\kern.37em\vrule height1.4ex width.05em
depth-.011em\kern-.37em C}}{\hbox{\sevenrm\kern.37em\vrule height1.3ex
width.05em depth-.02em\kern-.3em C}}{\hbox{\sevenrm\kern.37em\vrule height1.3ex
width.05em depth-.02em\kern-.3em C}}}}

\newcount\notitle
\notitle=1
\headline={\ifodd\pageno\rhead\else\lhead\fi}
\def\rhead{\ifnum\pageno=\notitle\hfill\else\hfill\eightrm\titolo\hfill
\folio\fi}
\def\lhead{\ifnum\pageno=\notitle\hfill\else\eightrm\folio\hfill\autore\hfill
\fi}
\newbox\bibliobox
\def\setref #1{\setbox\bibliobox=\hbox{[#1]\enspace}
        \parindent=\wd\bibliobox}
\def\biblap#1{\noindent\hang\rlap{[#1]\enspace}\indent\ignorespaces}
\def\art#1 #2: #3! #4! #5 #6 #7-#8 \par{\biblap{#1}#2: {\sl #3\/}.
        #4 {\bf #5}~(#6)\if.#7\else, \hbox{#7--#8}\fi.\par\smallskip}
\def\book#1 #2: #3! #4 \par{\biblap{#1}#2: {\bf #3.} #4.\par\smallskip}
\gdef\begin #1 #2\par{\xdef\titolo{#2}
\ifsect\let\neweq=\neweqt\else\let\neweq=\neweqf\fi
\ifsect\let\newthm=\newthmt\else\let\newthm=\newthmf\fi
\centerline{\titlefont\titolo}
\if\bititolo\empty\else\medskip\centerline{\titlefont\bititolo}
\xdef\titolo{\titolo\ \bititolo}\fi
\bigskip
\centerline{\bigfont by \autore}
\bigskip\bigskip
\ifsect\else\global\thmno=1\global\eqnumbo=1\fi}
\def\atanh{\mathop{\rm atanh}\nolimits}
\overfullrule=0pt
\let\smb=\smbar \let\j=\jmath 
\lineskiplimit=1pt

\def\autore{Marco Abate and Giorgio Patrizio}
\begin June Holomorphic curvature of Finsler metrics and complex geodesics
 
\bigskip
\smallsect 0. Introduction
 
If $D$ is a bounded convex domain in $\C^n$ , then the work of Lempert [L] and
Royden-Wong [RW] (see also [A]) show that given any point $p\in D$ and any
non-zero tangent vector $v\in\C^n$ at $p$, there exists a holomorphic map
$\phe\colon U\to D$ from the unit disk $U\subset\C$ into $D$ passing through
$p$ and tangent to $v$ in $p$ which is an isometry with respect to the
hyperbolic distance of $U$ and the Kobayashi distance of $D$. Furthermore if
$D$ is smooth and strongly convex then given $p$ and $v$ this holomorphic disk
is uniquely determined.
 
For a general complex manifold it is hard to determine whether or not such
complex curves, called complex geodesics in [V], exist. Therefore it is
natural to investigate the special properties enjoyed by the Kobayashi metric
of a strongly convex domain.  In this case it is known that the Kobayashi
metric is a strongly pseudoconvex smooth complex Finsler metric. Furthermore,
for a suitable notion of holomorphic curvature (see below), this metric in
convex domain has negative constant holomorphic curvature (cf.~[W], [S], [R]).
 
In [AP] it was started a systematic differential geometrical study of complex
geodesics in the framework of complex Finsler metrics. As in ordinary
Riemannian geometry it is natural to study geodesics as solutions of an
extremal problem and not as globally length-minimizing curves, so in our case
the natural notion turns out to be the one of {\sl geodesic complex curves,}
i.e., of holomorphic maps from the unit disk into the manifold sending
geodesics for the hyperbolic metric into geodesics for the given Finsler
metric. For instance, the annulus in $\C$ has no complex geodesics in the
sense of [V], whereas the usual universal covering map is a geodesic complex
curve in the previous sense.
 
It was shown in [AP] that geodesic complex curves for complex Finsler metrics
satisfy a system of partial differential equations and, under suitable
hypotheses, it was given an uniqueness theorem. Here we shall be concerned
with the question of existence.
 
The further ingredient needed to attack this problem is the notion of
holomorphic curvature of complex Finsler metrics. Given a complex manifold $M$
and a complex Finsler metric $F\colon T^{1,0}M\to\R$, i.e., a nonnegative
upper semicontinuous function such that
$$F(p;\lambda v) = |\lambda | F(p;v)$$
for all $(p;v)\in T^{1,0}M$ and $\lambda\in\C$, the holomorphic curvature of
$F$ at $p$ in the direction $v$ is the supremum of the Gaussian curvature at
the origin of the (pseudo)hermitian metrics on the unit disk obtained by
pulling back~$F$ via holomorphic maps $\phe\colon U\to M$ with $\phe(0)=p$ and
$\phe'(0)=\lambda v$ for some $\lambda\in\C^*$. Here, following [H], we make
the choice of computing the Gaussian curvature using the weak laplacian rather
then using Ahlfors' notion of supporting metrics; this approach seems to be
more natural for our applications and has a better connection with the usual
hermitian geometry. We remark that if $F$ is the norm associated to a
hermitian metric, then Wu in [Wu] showed that this definition yields the usual
holomorphic sectional curvature in the direction~$v$. We shall see in
section~2 that for smooth strongly pseudoconvex (see below for definitions)
Finsler metrics this definition recovers the
holomorphic sectional curvature defined by Kobayashi in~[K2].
 
In section~1 we give a survey of the elementary implications of this notion
of holomorphic curvature of complex Finsler metrics gotten by application of
Ahlfors' lemma and of its sharp form due to Heins [H]. For instance, as one
could expect, it is easy to prove a hyperbolicity criterion in terms of
negatively curved upper semicontinuous Finsler metric (cf. Corollary~1.5). As
for the investigation of geodesic complex curves and differential geometric
properties of intrinsic metrics, not much can be achieved without some
smoothness assumptions. One of the few general facts obtained here is the
explanation (Proposition~1.6) in terms of the holomorphic curvature of the
known property of the Carath\'eodory metric that a holomorphic disk which is
an isometry at one point is in fact an infinitesimal isometry at every point
(i.e., an infinitesimal complex geodesic in the terminology introduced in
[V]). Furthermore, in Proposition~1.7 it is given a very weak characterization
of the Kobayashi metric which nevertheless is useful later in the paper.
 
In order to prove more significant results it is necessary to consider the
case of smooth Finsler metric, i.e., such that $F\in
C^\infty\bigl(T^{1,0}M\setminus\{\hbox{Zero section}\}\bigr)$, which are in
addition strongly pseudoconvex, that is such that for every $p\in M$ the
indicatrix
$$I_p(M)=\{v\in T^{1,0}M\mid F(p;v) < 1\}$$
is strongly pseudoconvex. Under these assumptions in
section~2 we show how to compute the holomorphic curvature $K_F$ of $F$ by
means of a tensor explicitly defined in terms of $F$ and which agrees with the
usual one in case $F$ is the norm associated to a hermitian metric. The
expression we get has also been considered from a slightly different point of
view by Royden~[R] and Kobayashi~[K2].
 
After this preparation we address the problem of the existence of geodesic
complex curves in section~3. Let
$$A(\zeta)={-2\smbar\zeta\over1-|\zeta|^2}$$
for $\zeta\in U$, and
$$\Gamma^\alpha_{;i}=G^{\alpha\smbar\mu}G_{\smbar\mu;i},$$
where $G = F^2$, lower indices indicate derivatives (with respect to the
coordinates of the manifold those after semicolon, with respect to the
coordinates of the tangent space the others), $(G^{\alpha\smbar\beta}) =
(G_{\alpha\smbar\beta})^{-1}$and we are using the usual summation convention.
From results of [AP] it follows that the geodesic complex curves are
holomorphic maps $\phe\colon U\to M$ satisfying, for $\alpha= 1,\ldots, n$,
$$(\phe'')^\alpha+A\cdot(\phe')^\alpha+\Gamma^\alpha_{;i}(\phe;\phe')\phe^i=0,
        \neweq\eqzuno$$
and an additional set of equations which automatically hold if along the curve
the metric $F$ satisfies a K\"ahler condition introduced by Rund [Ru2], which
reduces to the usual one for hermitian metrics.
 
The holomorphic solutions of (0.1) have nice properties on their own. For
instance, they realize the holomorphic curvature at every point for the
direction tangent to the disk, and if they are isometry at the origin then
they are infinitesimal complex geodesics (cf. Proposition 3.2). Our first main
result (Theorem~3.3) describes necessary and sufficient conditions for the
holomorphic solvability of \eqzuno, and hence for the existence (and
uniqueness) of complex geodesics through a given point and direction.
 
The characterization given in Theorem~3.3 is completely expressed in terms of
the metric, but rather technical; to give a clearer geometric characterization
it is necessary to bring the curvature into the picture. The previous list of
properties of holomorphic solutions of \eqzuno\ shows that a natural necessary
condition is that $F$ has constant negative holomorphic curvature.
This is almost sufficient; to get the correct geometric conditions it is
necessary to introduce a further tensor, defined on the sphere bundle
$S^{1,0}M=\{(p;\xi)\in T^{1,0}M\mid F(p;\xi) = 1\}$. Set for
$\alpha= 1,\ldots, n$,
$$H_\alpha(v)=H_{\alpha i\smbar\mu\smbar\jmath}\,\bar{v^\mu}v^i\bar{v^j}=\bigl
        [(G_{\tau\smbar\mu}\Gamma^\tau_{;\alpha\smbar\jmath})_i-(G_{\tau\smbar
        \mu}\Gamma^\tau_{;i\smbar\jmath})_\alpha\bigr]\,\bar{v^\mu}v^i
        \bar{v^j},$$
where $\Gamma^\alpha_{;i\smbar\jmath}=(G^{\alpha\smbar\mu}G_{\smbar\mu;i})
_{;\smbar\jmath}$. To understand the meaning of this tensor, let us
consider the case of a hermitian metric. Then $H_{\alpha i\smb\mu\smbar\j}
=R_{i\alpha\smb\mu\smb\j}-R_{\alpha i\smb\mu\smb\j}$, where $R$ is the
Riemannian curvature tensor of the Chern connection associated to the
hermitian metric. In particular, $H_{\alpha i\smb\mu\smb\j}
\equiv0$ is equivalent to $\bar\de T\equiv 0$, where $T$ is the torsion form
of the Chern connection; $T$ is a
$T^{1,0}M$-valued $2$-form vanishing exactly when the given metric is
K\"ahler. In conclusion, $H$ may be interpreted as a torsion of the curvature;
in fact, it can be shown that $H_\alpha\equiv0$ is a simmetry condition ---
formally identical to (3.29) --- on the curvature tensor of the Chern connection
induced by the Finsler metric on the vertical subbundle of the two-fold
tangent bundle $T^{1,0}\bigl(T^{1,0}M\bigr)$. We intend to pursue these
matters elsewhere.
 
Using the tensor $H$ we may finally summarize our main results. We have
(Theorem~3.6):
 
\newthm Theorem \zuno: Let $M$ be a complex manifold  equipped with a strongly
pseudoconvex smooth complete Finsler metric $F$. Then there exists a unique
holomorphic solution $\phe\colon U\to M$ of $\eqzuno$ with $\phe(0) = p$ and
$\phe'(0) = \xi$ for any $(p;\xi)\in S^{1,0}M$ iff the holomorphic curvature
$K_F\equiv-4$ and $H_\alpha\equiv0$ for all $\alpha$.
 
In other words, the existence and uniqueness of holomorphic solutions of
\eqzuno\ is equivalent to constant negative holomorphic curvature and a
simmetry condition on a curvature tensor.
Furthermore, it turns out the the K\"ahler condition along a holomorphic
solution of \eqzuno\ holds iff it holds at one point, and thus (Corollary~3.9)
 
\newthm Theorem \zdue: Let $M$ be a complex manifold equipped with a strongly
pseudoconvex smooth complete Finsler metric $F$. Assume that the holomorphic
curvature $K_F\equiv-4$ and that $H_\alpha\equiv0$ for all $\alpha$. Take
$(p_0;\xi_0)\in S^{1,0}M$. Then there is a (a fortiori unique) geodesic
complex curve passing through $p_0$ tangent to $\xi_0$ iff the K\"ahler
condition holds at~$(p_0;\xi_0)$.
 
From this result it is also possible to obtain a geometric characterization of
the Kobayashi metric. The following corollary (Corollary~3.10) improves
results of Pang [P] and it is closely related to those of Faran [F]:
 
\newthm Corollary \ztre: Let $F$ be a strongly pseudoconvex smooth complete
Finsler metric with constant holomorphic curvature $K_F\equiv-4$ and such that
$H_\alpha\equiv0$ for all $\alpha$. Then $F$ is the Kobayashi metric of $M$.
 
The second named author thanks the Max-Planck-Institut f\"ur Mathematik of
Bonn for its hospitality and support while this paper was completed.
 
\smallsect 1. Holomorphic curvature for semicontinuous metrics
 
Let $U$ denote the unit disk in the complex plane. A {\sl pseudohermitian
metric} $\mu_g$ of {\sl scale}~$g$ on~$U$ is the upper semicontinuous
pseudometric on the tangent bundle of~$U$ defined by
$$\mu_g=g\,d\zeta\otimes d\smbar\zeta,\neweq\equuno$$
where $g\colon U\to\R^+$ is a non-negative upper semicontinuous function such
that $S_g=g^{-1}(0)$ is a discrete subset of~$U$.
 
If $\mu_g$ is a standard hermitian metric on~$U$, i.e., $g$ is a $C^2$
positive function, the Gaussian curvature of~$\mu_g$ is defined by
$$K(\mu_g)=-{1\over 2g}\Delta\log g,\neweq\equdue$$
where $\Delta$ denotes the usual Laplacian
$$\Delta u=4\,{\de^2 u\over\de\zeta\de\smbar\zeta}.\neweq\equtre$$
The ({\sl lower}) {\sl generalized Laplacian} of an upper semicontinuous
function~$u$ is defined by
$$\Delta u(\zeta)=4\liminf_{r\to0}{1\over r^2}\left\{{1\over2\pi}\int_0^{2\pi}
        u(\zeta+re^{i\theta})\,d\theta-u(\zeta)\right\}.\neweq\eququattro$$
It is worthwhile to remark explicitely some features of this definition. First
of all, when $u$ is a function of class~$C^2$ in a neighbourhood of the
point~$\zeta_0$, \eququattro\ actually reduces to~\equtre. In fact, for $r$
small enough we can write
$$\eqalign{u(\zeta_0+re^{i\theta})=u(\zeta_0)&+{\de u\over\de\zeta}(\zeta_0)r
        e^{i\theta}+{\de u\over\de\smbar\zeta}(\zeta_0)re^{-i\theta}\cr
        &+{1\over2}{\de^2 u\over\de\zeta^2}(\zeta_0)r^2e^{2i\theta}+
        {\de^2u\over\de\zeta\de\bar\zeta}(\zeta_0)r^2+{1\over2}{\de^2u\over
        \de\bar\zeta^2}(\zeta_0)r^2e^{-2i\theta}+o(r^2);\cr}$$
hence
$${1\over2\pi}\int_0^{2\pi}u(\zeta_0+re^{i\theta})\,d\theta-u(\zeta)=r^2{\de^2
        u\over\de\zeta\de\smbar\zeta}(\zeta_0)+o(r^2),$$
and the claim follows.
 
Second, if $u$ is an upper semicontinuous function, it is easy to see that
$\Delta u\ge0$ is equivalent to the submean property; so $\Delta u\ge0$ iff
$u$ is subharmonic.
 
Finally, if $\zeta_0$ is a point of local maximum for $u$, then clearly
$\Delta u(\zeta_0)\le0$.
 
Now let $\mu_g$ be a pseudohermitian metric on $U$. Then the {\sl Gaussian
curvature} $K(\mu_g)$ of~$\mu_g$ is the function defined on $U\setminus S_g$
by \equdue\ --- using the generalized Laplacian \eququattro; clearly, if $\mu_
g$ is a standard hermitian metric, $K(\mu_g)$ reduced to the usual Gaussian
curvature.
 
The idea behind the classical Ahlfors lemma is to
compare a generic pseudohermitian metric with an extremal one --- usually the
Poincar\'e metric. For $a>0$, let $g_a\colon U\to\R^+$ be defined by
$$g_a(\zeta)={1\over a(1-|\zeta|^2)^2};$$
then $\mu_a=g_a\,d\zeta\otimes d\smbar\zeta$ is a hermitian metric of
constant Gaussian curvature $K(\mu_a)=-4a$. Of course, $\mu_1$ is the standard
Poincar\'e metric on~$U$.
 
The classical Ahlfors lemma is true in this more general situation:
 
\newthm Proposition \uuno: Let $\mu_g=g\,d\zeta\otimes d\smbar\zeta$ be a
pseudohermitian metric on~$U$ such that $K(\mu_g)\le-4a$ on $U\setminus S_g$
for some $a>0$. Then $g\le g_a$.
 
\pf The proof follows closely the classical one due to Ahlfors. For the sake
of completeness we report it here.
 
For $0<r<1$, define $U_r=\{\zeta\in\C\mid|\zeta|<r\}$ and $g_a^r\colon U_r\to
\R^+$ by
$$g_a^r(\zeta)={1\over a(1-|\zeta|^2/r^2)^2}=g_a(\zeta/r),$$
and set $f_r=g/g_a^r\colon U_r\to\R^+$. Being upper semicontinuous, $g$ is
bounded on~$\bar{U_r}$; since $g_a^r(\zeta)\to+\infty$ as $|\zeta|\to r$, there
is a point $\zeta_0\in U_r$ of maximum for~$f_r$. Clearly, $\zeta_0\notin S_g$;
hence
$$0\ge\Delta\log f_r(\zeta_0)\ge\Delta\log g(\zeta_0)-\Delta\log g_a^r(\zeta_0)
        =-2g(\zeta_0)K(\mu_g)(\zeta_0)-8ag_a^r(\zeta_0).\neweq\equsette$$
By assumption, $K(\mu_g)(\zeta_0)\le-4a$; therefore \equsette\ yields $g(\zeta_
0)\le g_a^r(\zeta_0)$ and thus, by the choice of~$\zeta_0$,
$$\forevery{\zeta\in U_r} g(\zeta)\le g_a(\zeta/r).$$
Letting $r\to 1^-$ we obtain the assertion.\qedn
 
To complement this result, we recall a theorem due to Heins [H, Theorem~7.1],
showing that $\mu_g\ne\mu_a$ in the statement of Proposition~\uuno\ implies
that $g$ is strictly less than~$g_a$ everywhere:
 
\newthm Theorem \udue: (Heins) Let $\mu_g=g\,d\zeta\otimes d\smbar\zeta$ be a
pseudohermitian metric on~$U$ such that $K(\mu_g)\le-4a$ on~$U\setminus S_g$
for some~$a>0$. Assume there is $\zeta_0\in U\setminus S_g$ such that
$g(\zeta_0)=g_a(\zeta_0)$. Then $\mu_g\equiv\mu_a$.
 
Now we start looking to the several variables situation. If $M$ is a complex
manifold, we shall denote by $TM$ its real tangent bundle endowed
with the almost-complex structure~$J$ induced by the complex structure of~$M$;
by $T^cM$ the complexification of $TM$ and by
$T^{1,0}M$ the $(1,0)$-part (i.e., the $i$-eigenspace of~$J$) of $T^cM$.
As well known, $T^{1,0}M$ is naturally complex-isomorphic to $TM$. In this
paper we shall mainly use $T^{1,0}M$ as representative of the tangent bundle
of~$M$, except for an argument needed in section~3.
 
A {\sl complex Finsler metric}~$F$ on a complex $n$-dimensional ($n\ge1$)
manifold~$M$ is an upper semicontinuous map $F\colon T^{1,0}M\to\R^+$
satisfying
\smallskip
\item{(i)}$F(p;v)>0$ for all $p\in M$ and $v\in T^{1,0}_pM$ with $v\ne0$;
\item{(ii)}$F(p;\lambda v)=|\lambda|F(p;v)$ for all $p\in M$, $v\in T^{1,0}_pM$
and $\lambda\in\C$.
\smallskip
\noindent We shall sistematically denote by $G\colon T^{1,0}(M)\to\R^+$ the
function $G=F^2$. Note that, thanks to condition (ii), the definition of
length of a smooth curve in a Riemannian manifold makes sense in this context
too; so we may again associate to $F$ a topological distance on~$M$, and we
shall say that $F$ is {\sl complete\/} if this distance is. For the same
reason, it makes sense to call {\sl (real) geodesics} the extremals of the
length functional. General introductions to real Finsler geometry are~[Ru1, B].
 
Take $p\in M$ and $v\in T^{1,0}_pM$, $v\ne0$. The {\sl
holomorphic curvature} $K_F(p;v)$ of~$F$ at~$(p;v)$ is given by
$$K_F(p;v)=\sup\{K(\phe^*G)(0)\},$$
where the supremum is taken with respect to the family of all holomorphic maps
$\phe\colon U\to M$ with $\phe(0)=p$ and $\phe'(0)=\lambda v$ for some
$\lambda\in\C^*$, and $K(\phe^*G)$ is the Gaussian curvature discussed so
far of the pseudohermitian metric $\phe^*G$ on~$U$. Clearly, the holomorphic
curvature depends only on the complex line in~$T^{1,0}_pM$ spanned by~$v$, and
not on~$v$ itself.
 
The holomorphic curvature may also be defined  (see e.g. [S]) taking the
supremum with respect to the family of all holomorphic maps $\phe\colon U_r\to
M$ with $\phe(0)=p$ and $\phe'(0)=v$, where $U_r\subset\C$ is the disk of
center the origin and radius~$r$. We chose the given definition to stress the
similarities with the definitions of the Kobayashi and Carath\'eodory metrics.
 
If $F$ is a complex Finsler metric on~$U$ --- and so $G=F^2$ is a
pseudohermitian metric $G= g\,d\zeta\otimes d\smbar\zeta$ on~$U$ ---, a
priori we have defined two curvatures for~$F$: $K_F(\zeta;1)$ and
$K(G)(\zeta)$. As anybody may guess, they actually coincide; this is a
consequence of
 
\newthm Lemma \utre: Let $\mu_g=g\,d\zeta\otimes d\smbar\zeta$ be a
pseudohermitian metric on~$U$, and $\phe\colon U\to U$ a holomorphic self-map
of~$U$. Then on $U\setminus[\phe^{-1}(S_g)\cup S_{\phe'}]$
$$K(\phe^*\mu_g)=K(\mu_g)\circ\phe.$$
 
\pf Let $\{g_n\}$ be a sequence of $C^2$ functions such that $g_n\ge g_{n+1}$
and with $g_n(x)\to g(x)$. Then on $U\setminus[\phe^{-1}(S_{g_n})\cup S_{\phe'}]
\supset U\setminus[\phe^{-1}(S_g)\cup S_{\phe'}]$ we have
$$\eqalign{K(\phe^*\mu_{g_n})&=-{1\over2|\phe'|^2(g_n\circ\phe)}\Delta\log
        (|\phe'|^2g_n\circ\phe)\cr
        &=-{2\over|\phe'|^2(g_n\circ\phe)}\left[{\de^2\log(g_n\circ\phe)\over
        \de\zeta\de\smbar\zeta}+{\de^2\log|\phe'|^2\over\de\zeta\de\smbar\zeta}
        \right]\cr
        &=-{2\over|\phe'|^2(g_n\circ\phe)}\left[|\phe'|^2\left({\de^2\log g_n
        \over\de\zeta\de\smbar\zeta}\right)\circ\phe\right]\cr
        &=-{1\over2(g_n\circ\phe)}\bigl(\Delta\log g_n\bigr)\circ\phe=
        K(\mu_{g_n})\circ\phe.\cr}$$
Letting $n\to+\infty$ and applying the dominated convergence theorem we get
the assertion.\qedn
 
The holomorphic curvature defined in this way is clearly invariant under
holomorphic isometries. More generally, if $f\colon M\to N$ is holomorphic and
$F$ is a Finsler metric on~$N$, we have
$$K_F\bigl(f(z);df_z(v)\bigr)\ge K_{f^*F}(z;v),$$
that is $f^*K_F\ge K_{f^*F}$.
 
When $F$ is a honest smooth hermitian metric on~$M$, $K_F(p;v)$ coincides with
the usual holomorphic sectional curvature of~$F$ at~$(p;v)$ (see~[Wu]). The
aim of this section is to extend a couple of results already known for
hermitian metrics to this more general case.
 
A piece of terminology: we say that a complex Finsler metric $F$ has {\sl
holomorphic curvature bounded above (below) by a constant} $c\in\R$ if $K_
F(p;v)\le c$ (respectively, $K_F(p;v)\ge c$) for all $(p;v)\in T^{1,0}M$ with
$v\ne0$.
 
Our first result is the usual several variables version of Ahlfors' lemma:
 
\newthm Proposition \duno: Let $F$ be a complex Finsler metric on a complex
manifold~$M$. Assume that the holomorphic curvature of~$F$ is bounded above by
a negative constant~$-4a$, for some~$a>0$. Then
$$\phe^*F\le\mu_a\neweq\eqddue$$
for all holomorphic maps $\phe\colon U\to M$.
 
\pf $\phe^*F$ is a pseudohermitian metric on~$U$; by assumption (and by
Lemma~\utre), $K(\phe^*F)\le-4a$. Then the assertion follows from
Proposition~\uuno.\qedn
 
As a consequence, we obtain a generalization of a well-known criterion of
hyperbolicity:
 
\newthm Corollary \ddue: Let $M$ be a complex manifold admitting a (complete)
complex Finsler metric~$F$ with holomorphic curvature bounded above by a
negative constant. Then $M$ is (complete) hyperbolic.
 
\pf Up to multiplying $F$ by a suitable constant, we may assume $K_F\le-4$.
Let $d$ denote the distance induced by~$F$ on~$M$, and $\omega$ the
Poincar\'e distance on~$U$. Then Proposition~\duno\ yields
$$d\bigl(\phe(\zeta_1),\phe(\zeta_2)\bigr)\le\omega(\zeta_1,\zeta_2),$$
for all $\zeta_1$,~$\zeta_2\in U$ and holomorphic maps $\phe\colon U\to M$.
But this immediately implies (cf.~[K1,~Proposition~I\negthinspace V.1.4]) that
the Kobayashi distance~$k_M$ of~$M$ is bounded below by~$d$, and the assertion
follows.\qedn
 
In particular, then, a complex manifold admitting a complete complex Finsler
metric with holomorphic curvature bounded above by a negative constant is
necessarily taut.
 
The notion of holomorphic curvature for (non-smooth) Finsler metric has been
introduced recently in connection with the Carath\'eodory and Kobayashi
metrics. In particular, Wong [W] and Suzuki [S] (see also [Bu]) have shown
that the holomorphic curvature of the Carath\'eodory metric is bounded {\sl
above} by~$ -4$, whereas the holomorphic curvature of the Kobayashi metric is
bounded {\sl below} by~$-4$.
 
An interesting immediate consequence of this is an interpretation in terms of
curvature of a well known property of the Carath\'eodory metric:
 
\newthm Proposition \usei: Let $F$ be a complex Finsler metric on a manifold
$M$ with holomorphic curvature bounded above by~$-4$. Let $\phe\colon U\to M$
be a holomorphic map. Then the following are equivalent:
{\smallskip
\itm{(i)}$\phe^*F(0;1)=F\bigl(\phe(0);\phe'(0)\bigr)=1$, that is $\phe$ is an
isometry at the origin between the Poincar\'e metric on $U$ and~$F$;
\itm{(ii)}$\phe$ is an infinitesimal complex geodesic, that is $\phe^*F$ is
the Poincar\'e metric of~$U$.}
 
\pf By definition and Lemma~\utre, the Gaussian curvature of $\phe^*F$ is
bounded above by~$-4$. The assertion follows from Heins' Theorem~\udue.\qedn
 
Bounds on the holomorphic curvature allow to compare a complex Finsler metric
to the Kobayashi metric --- and maybe to prove that a given Finsler metric
actually is the Kobayashi metric. For instance, Pang [P] and Faran [F] gave
conditions under which a smooth complex Finsler metric of constant negative
holomorphic curvature coincides with the Kobayashi metric. We shall discuss
the smooth case in detail in the next two sections; here, to provide the right
set-up to the problem, we examine a bit the general situation.
 
We need an auxiliary notion to formulate our observation. Let $F$ be a complex
Finsler metric on a manifold~$M$, and take $(p;v)\in T^{1,0}M$. We
say that $F$ is {\sl realizable} at~$(p;v)$ if there is a holomorphic map
$\phe\colon U\to M$ such that $\phe(0)=p$ and $\lambda\phe'(0)=v$ with
$|\lambda|=F(p;v)$. In other words, $\phe$ is an isometry at the origin between
the Poincar\'e metric of~$U$ and~$F$.
 
Obviously, the Kobayashi metric is realizable in any taut manifold; on the
other hand, as a consequence of the next result, the Carath\'eodory metric is
realizable iff it coincides with the Kobayashi metric.
 
\newthm Proposition \usette: Let $F$ be a complex Finsler metric on a
manifold~$M$, and choose $(p_0;v_0)\in T^{1,0}M$. Then:
{\smallskip
\itm{(i)}If $F$ is realizable at $(p_0;v_0)$, then $F(p_0;v_0)\ge\kappa_M(p_0;
v_0)$;
\itm{(ii)}If $K_F\le-4$, then $F\le\kappa_M$;
\itm{(iii)}If $F$ is realizable at $(p_0;v_0)$ and $K_F\le-4$, then
$F(p_0;v_0)=\kappa_M(p_0;v_0)$.}
 
\pf (i) Let $\phe\colon U\to M$ be a holomorphic map with $\phe(0)=p_0$ and
$v_0=\lambda\phe'(0)$ such that $F(p_0;v_0)=|\lambda|$. Then
$$F(p_0;v_0)=|\lambda|\ge\kappa_M(p_0;v_0).$$
 
(ii) Take $(p;v)\in T^{1,0}M$ and let $\phe\colon U\to M$ be a holomorphic map
with $\phe(0)=p$ and $v=\lambda\phe'(0)$. Then $\phe^*G$ is a pseudohermitian
metric on~$U$ with
Gaussian curvature bounded above by~$-4$ (by Lemma~\utre); it follows from
Proposition~\uuno\ that $\phe^*G\le\mu_1$. Thus
$$F(p;v)=F\bigl(\phe(0);\lambda\phe'(0)\bigr)=\phe^*F(0;\lambda)\le|\lambda|.$$
Since this holds for all such $\phe$, we get $F\le\kappa_M$.
 
(iii) Obvious, now.\qedn
 
This is the best it can be done on the basis of Ahlfors' lemma and Heins'
theorem only. In order to get deeper results it is necessary to use more tools
--- as we shall see in the smooth case discussed in the rest of the paper.
 
\smallsect 2. Holomorphic curvature: the smooth case
 
In this section we shall derive a tensor expression of the holomorphic
curvature of a smooth complex Finsler metric.
 
First of all, we need a few definitions, notations and general formulas. Let
$F$ be a complex Finsler metric on a complex manifold $M$, and set $G=F^2$, as
usual. We shall assume that $F$ is {\sl smooth,} that is that
$F$ is of class $C^k$ ($k\ge4$) out of the zero section of $T^{1,0}M$. By
the way, $(T^{1,0}M)_0$ will denote the complement in $T^{1,0}M$ of the zero
section.
 
If $(z^1,\ldots,z^n)$ are local coordinates on $M$, a local section
of $T^{1,0}M$ will be written as
$$\sum_{j=1}^n v^j{\de\over\de z^j},$$
and we shall use $(z^1,\ldots,z^n;v^1,\ldots,v^n)$ as local coordinates on $T^
{1,0}M$.
 
We shall denote by indexes like $\alpha$, $\smb\beta$ and so on the
derivatives with respect to the $v$-coordinates; for instance,
$$G_{\alpha\smb\beta}={\de^2 G\over\de v^\alpha\de\bar{v^\beta}}.$$
On the other hand, the derivatives with respect to the $z$-coordinates will be
denoted by indexes after a semicolon; for instance,
$$G_{;ij}={\de^2 G\over\de z^i\de z^j}\qquad\hbox{or}\qquad G_{\alpha;\smb\j}
        ={\de^2 G\over\de\bar{z^j}\de v^\alpha}.$$
 
A smooth complex Finsler metric $F$ will be said {\sl strongly pseudoconvex} if
the $F$-indicatrices are strongly pseudoconvexes, i.e., if the Levi matrix
$(G_{\alpha\smb\beta})$ is positive definite on $(T^{1,0}M)_0$. As usual in
hermitian geometry, we shall denote by $(G^{\alpha\smb\beta})$ the inverse
matrix of $(G_{\alpha\smb\beta})$, and we shall use it to raise indexes.
The usual summation convention will hold throughout the rest of the paper.
 
The main (actually, almost the unique) property of the function $G$ is its
(1,1)-homo\-geneity: we have
$$G(z;\lambda v)=\lambda\smb\lambda\,G(z;v)\neweq\eqdhom$$
for all $(z;v)\in T^{1,0}M$ and $\lambda\in\C$. We now collect a number of
formulas we shall use later on which follows from \eqdhom. First of
all, differentiating with respect to $v^\alpha$ and $\bar{v^\beta}$ we get
$$\eqalign{G_\alpha(z;\lambda v)&=\smb\lambda G_\alpha(z;v),\cr
        G_{\alpha\smb\beta}(z;\lambda v)&=G_{\alpha\smb\beta}(z;v),\cr
        G_{\alpha\beta}(z;\lambda v)&=(\smb\lambda/\lambda)G_{\alpha\beta}
        (z;v).\cr}$$
Thus differentiating with respect to $\lambda$ or $\smb\lambda$ and then
setting $\lambda=1$ we get
$$G_{\alpha\smb\beta}\,\bar{v^\beta}=G_\alpha,\qquad
        G_{\alpha\beta}\,v^\beta=0,\neweq\eqdf$$
and
$$G_{\alpha\beta\gamma}\,v^\gamma=-G_{\alpha\beta},\qquad
        G_{\alpha\beta\smb\gamma}\,\bar{v^\gamma}=G_{\alpha\beta},\qquad
        G_{\alpha\smb\beta\gamma}\,v^\gamma=0,\neweq\eqde$$
where everything is evaluated at $(z;v)$.
 
On the other hand, differentiating directly \eqdhom\ with respect to $\lambda$
or $\smb\lambda$ and putting eventually $\lambda=1$ we get
$$G_\alpha\,v^\alpha=G,\qquad G_{\alpha\beta}\,v^\alpha v^\beta=0,\qquad
        G_{\alpha\smb\beta}\,v^\alpha\bar{v^\beta}=G.\neweq\eqdb$$
It is clear that we may get other formulas applying any differential operator
acting only on the $z$-coordinates, or just by conjugation. For instance, we
get
$$G_{\smb\alpha;i}\,\bar{v^\alpha}=G_{;i},\neweq\eqdo$$
and so on.
 
Assuming now $F$ strongly pseudoconvex, we get another bunch of formulas we
shall need later on. First of all, applying $G^{\alpha\smb\beta}$ to the first
equation in \eqdf\ we get
$$G^{\alpha\smb\beta}G_\alpha=\bar{v^\beta},\neweq\eqda$$
and thus, applying \eqdo,
$$G_{\smb\beta;i}G^{\alpha\smb\beta}G_{\alpha}=G_{;i}.\neweq\eqdd$$
Recalling that $(G^{\alpha\smb\beta})$ is the inverse matrix of $(G_{\alpha
\smb\beta})$, we may also compute derivatives of~$G^{\alpha\smb\beta}$:
$$DG^{\alpha\smb\beta}=-G^{\alpha\smb\nu}G^{\mu\smb\beta}(DG_{\mu\smb\nu}),
        \neweq\eqdc$$
where $D$ denotes any first order linear differential operator. As a
consequence of \eqde\ and~\eqdc\ we get
$$G^{\alpha\smb\beta}_{\smb\sigma}\,\bar{v^\sigma}=-G^{\alpha\smb\nu}G^{\mu
        \smb\beta}G_{\mu\smb\nu\smb\sigma}\,\bar{v^\sigma}=0,\neweq\eqdg$$
and recalling also \eqda\ we obtain
$$G_{\smb\beta}G^{\alpha\smb\beta}_\gamma=-G_{\smb\beta}G^{\mu\smb\beta}
        G^{\alpha\smb\nu}G_{\mu\smb\nu\gamma}=-G^{\alpha\smb\nu}G_{\mu\smb\nu
        \gamma}v^\mu=0.\neweq\eqdh$$
 
Now we may start to work. Our first goal is to compute the holomorphic
curvature of our strongly pseudoconvex smooth complex Finsler metric $F$. Set
$$S^{1,0}M=\{\xi\in T^{1,0}M\mid F(\xi)=1\},$$
and choose $p\in M$ and $\xi\in S^{1,0}_pM$. To compute $K_F(p;\xi)$ we
should write the Gaussian curvature at the origin of $\phe^*G$, where
$\phe\colon U\to M$ is any holomorphic map with $\phe(0)=p$ and
$v=\phe'(0)=\lambda\xi$, where $|\lambda|=F(p;v)=[\phe^*G(0;1)]^{1/2}$.
 
Writing $\phe^*G=g\,d\zeta\otimes d\smb\zeta$, we have
$$g(\zeta)=G\bigl(\phe(\zeta);\phe'(\zeta)\bigr),\qquad g(0)=|\lambda|^2,$$
and
$$K(\phe^*G)(0)=-{1\over2g(0)}(\Delta\log g)(0)=-{2\over|\lambda|^2}
        {\de^2(\log g)\over\de\smb\zeta\de\zeta}(0).$$
The computation of the Laplacian yields
$$\eqalign{{\de^2(\log g)\over\de\smb\zeta\de\zeta}=&-{1\over G(\phe;\phe')^2}
        \Bigl|G_{;i}(\phe;\phe')(\phe')^i+G_\alpha(\phe;\phe')(\phe'')^\alpha
        \Bigr|^2\cr
        &+{1\over G(\phe;\phe')}\Bigl\{G_{;i\smb\j}(\phe;\phe')(\phe')^i
        \bar{(\phe')^j}+G_{\alpha\smb\beta}(\phe;\phe')(\phe'')^\alpha
        \bar{(\phe'')^\beta}\cr
        &\qquad\qquad\qquad\qquad+2\Re\bigl[G_{\smb\alpha;i}(\phe;\phe')
        (\phe')^i\bar{(\phe'')^\alpha}\bigr]\Bigr\}.\cr}$$
Hence writing $\eta=\phe''(0)$ we get
$$\eqalign{K(\phe^*G)(0)=&-2\bigl[G_{;i\smb\j}(p;\xi)-G_{;i}(p;\xi)G_{;\smb\j}
        (p;\xi)\bigr]\xi^i\bar{\xi^j}\cr
        &-{2\over|\lambda|^4}\bigl[G_{\alpha\smb\beta}(p;\xi)-G_\alpha(p;\xi)
        G_{\smb\beta}(p;\xi)\bigr]\eta^\alpha\bar{\eta^\beta}\cr
        &-{4\over|\lambda|^4}\Re\bigl\{\lambda^2\bigl[G_{\smb\alpha;i}(p;\xi)
        -G_{\smb\alpha}(p;\xi)G_{;i}(p;\xi)\bigr]\xi^i\bar{\eta^\alpha}\bigr\}
        .\cr}\neweq\eqdzero$$
We must compute the supremum (with respect to $\lambda$ and $\eta$) of this
formula. For the moment, let us consider $\lambda$ fixed, and look for the
infimum of
$$\bigl[G_{\alpha\smb\beta}(p;\xi)-G_\alpha(p;\xi)G_{\smb\beta}(p;\xi)\bigr]
        \eta^\alpha\bar{\eta^\beta}+2\Re\bigl\{\lambda^2\bigl[G_{\smb\alpha;i}
        (p;\xi)-G_{\smb\alpha}(p;\xi)G_{;i}(p;\xi)\bigr]\xi^i
        \bar{\eta^\alpha}\bigr\},$$
that is of
$$I_\lambda(\eta)=A_{\alpha\smb\beta}\,\eta^\alpha\bar{\eta^\beta}+2\Re\bigl[
        \lambda^2B_{\smb\alpha;i}\,\xi^i\bar{\eta^\alpha}\bigr],\neweq\eqduno$$
where
$$\displaylines{A_{\alpha\smb\beta}=G_{\alpha\smb\beta}(p;\xi)-G_\alpha(p;\xi)
        G_{\smb\beta}(p;\xi),\cr
        B_{\smb\alpha;i}=G_{\smb\alpha;i}(p;\xi)-G_{\smb\alpha}(p;\xi)
        G_{;i}(p;\xi).\cr}$$
Let us study the hermitian form $(A_{\alpha\smb\beta})$. By assumption, the
matrix $\bigl(G_{\alpha\smb\beta}(p;\xi)\bigr)$ induces a positive definite
hermitian product on $\C^n$; so we may decompose $\C^n$ accordingly as the
orthogonal sum of $\C\xi$ and its orthogonal $(\C\xi)^\perp$. Since, by \eqdb\
and \eqdo,
$$\displaylines{A_{\alpha\smb\beta}\,\eta^\alpha\bar{\xi^\beta}=G_{\alpha\smb
        \beta}\,\eta^\alpha\bar{\xi^\beta}-G_\alpha G_{\smb\beta}\,\eta^\alpha
        \bar{\xi^\beta}=G_\alpha\eta^\alpha-G_\alpha\eta^\alpha=0,\cr
        B_{\smb\alpha;i}\,\xi^i\bar{\xi^\alpha}=G_{\smb\alpha;i}\,\xi^i
        \bar{\xi^\alpha}-G_{\smb\alpha}G_{;i}\,\xi^i\bar{\xi^\alpha}=G_{;i}\xi^i
        -G_{;i}\xi^i=0,\cr}$$
for every $\eta\in\C^n$, if we denote by
$$\tilde\eta=\eta-\left(G_{\alpha\smb\beta}\,\eta^\alpha\bar{\xi^\beta}\right)
        \xi$$
the orthogonal projection of $\eta$ into $(\C\xi)^\perp$, we get
 
\newthm Lemma \duno: Let $F$ be a strongly pseudoconvex smooth complex Finsler
metric on $M$, and take $p\in M$ and $\xi\in S^{1,0}_pM$. Then
$I_{\lambda}\equiv 0$ on $\C\xi$ and $I_\lambda(\tilde\eta)=I_\lambda(\eta)$
for all $\eta\in\C^n$.
 
So it suffices to study $I_\lambda$ on $(\C\xi)^\perp$. Note that $\tilde\eta
\in(\C\xi)^\perp$ iff
$$0=G_{\alpha\smb\beta}\,\tilde\eta^\alpha\bar{\xi^\beta}=G_\alpha\tilde
        \eta^\alpha;$$
therefore on $(\C\xi)^\perp$ we have
$$A_{\alpha\smb\beta}\,\tilde\eta^\alpha=G_{\alpha\smb\beta}\,\tilde\eta^\alpha.
        \neweq\eqddue$$
In particular $(A_{\alpha\smb\beta})$ is positive definite on $(\C\xi)^\perp$.
 
Thus $I_\lambda$ is a quadratic polynomial on $(\C\xi)^\perp$ with positive
definite leading term; hence $I_\lambda$ attains a minimum at $\tilde\eta\in
(\C\xi)^\perp$ given by
$$A_{\alpha\smb\beta}\,\tilde\eta^\alpha=-\lambda^2B_{\smb\beta;i}\,\xi^i,\qquad
        \beta=1,\ldots,n,$$
that is, by \eqddue,
$$\tilde\eta^\alpha=-\lambda^2G^{\alpha\smb\beta}B_{\smb\beta;i}\,\xi^i,\qquad
        \alpha=1,\ldots,n.\neweq\eqdtre$$
Putting \eqdtre\ into \eqduno\ we find that the minimum of $I_\lambda$ is
$$-|\lambda|^4G^{\alpha\smb\beta}\bigl(B_{\alpha;\smb\j}\,\bar{\xi^j}\bigr)
        \bigl(B_{\smb\beta;i}\,\xi^i\bigr)<0,$$
and thus \eqdzero\ yields
$$K_F(p;\xi)=-2\bigl[G_{;i\smb\j}-G_{;i}G_{;\smb\j}-G^{\alpha\smb\beta}
        B_{\alpha;\smb\j}B_{\smb\beta;i}\bigr]\xi^i\bar{\xi^j}.$$
Now, using \eqdb\ and \eqda\ we get
$$G^{\alpha\smb\beta}B_{\alpha;\smb\j}B_{\smb\beta;i}=G^{\alpha\smb\beta}
        G_{\alpha;\smb\j}G_{\smb\beta;i}-G_{;i}G_{;\smb\j};$$
therefore
$$K_F(p;\xi)=-2\bigl[G_{;i\smb\j}-G^{\alpha\smb\beta}G_{\alpha;\smb\j}
        G_{\smb\beta;i}\bigr]\xi^i\bar{\xi^j}.\neweq\eqdquattro$$
It is easy to check (cf. [Wu]) that when $F$ is a standard hermitian metric on
$M$, then \eqdquattro\ reduces to the usual holomorphic sectional curvature in
the direction~$\xi$. Furthermore, \eqdquattro\ exactly yields the holomorphic
sectional curvature introduced by Kobayashi in~[K2].
 
There is a shorter way of writing $K_F$. Set
$$\Gamma^\alpha_{;i}=G^{\alpha\smb\mu}G_{\smb\mu;i},$$
and put $\Gamma^\alpha_{;i\smb\j}=(\Gamma^
\alpha_{;i})_{;\smb\j}$; then, by \eqdc,
$$\Gamma^\alpha_{;i\smb\j}=G^{\alpha\smb\mu}G_{\smb\mu;i\smb\j}-G^{\alpha
        \smb\nu}G^{\beta\smb\mu}G_{\beta\smb\nu;\smb\j}G_{\smb\mu;i},$$
and so, by \eqda\ and \eqdd,
$$G_\alpha\Gamma^\alpha_{;i\smb\j}=G_{;i\smb\j}-G^{\beta\smb\mu}G_{\beta;\smb\j}
        G_{\smb\mu;i}.\neweq\eqduseibis$$
Summing up, we have proved the
 
\newthm Proposition \ddue: Let $F$ be a strongly pseudoconvex smooth complex
Finsler metric on $M$, and take $(p;\xi)\in S^{1,0}M$. Then the holomorphic
curvature of $F$ in the direction of $\xi$ is
$$K_F(p;\xi)=-2G_\alpha\Gamma^\alpha_{;i\smb\j}\,\xi^i\bar{\xi^j}.
        \neweq\eqdcinque$$
 
For future reference, we note here that more generally the holomorphic
curvature of $F$ in the direction of a non-zero vector $v\in T^{1,0}_pM$
--- which coincides with the holomorphic curvature in the direction of
$\xi=v/F(p;v)$ --- is given by the formula
$$K_F(p;v)=-{2\over G(p;v)^2}G_\alpha(p;v)\Gamma^\alpha_{;i\smb\j}(p;v)\,
        v^i\bar{v^j}.\neweq\eqdcinquemez$$
 
\smallsect 3. Holomorphic curvature and geodesic complex curves
 
Let $F$ be a (smooth) complex Finsler metric on a manifold~$M$. Let $U_r$
denote the disk $\{\zeta\in\C\mid|\zeta|<r\}$ in $\C$ (with $0<r\le 1$),
endowed with the restriction of the Poincar\'e metric of $U$; note that $U_r$
is a convex subset of~$U$ with respect to the Poincar\'e metric. A holomorphic
map $\phe\colon U_r\to M$ is a {\sl segment of infinitesimal complex geodesic}
if $\phe^*F$ is the Poincar\'e metric on~$U_r$, that is if $\phe$ is a local
isometry from the Poincar\'e metric to~$F$. On the other hand, $\phe$ is said
{\sl segment of geodesic complex curve} if the image via $\phe$ of any (real)
geodesic in $U_r$ is a (real) geodesic for $F$ in $M$. In other words, $\phe$
is a local isometry and $\phe(U_r)$ is a totally geodesic complex curve in
$M$. When $r=1$, we shall talk of infinitesimal complex geodesics and geodesic
complex curves {\it tout-court.} In any case, if $(p;v)
=\bigl(\phe(0);\phe'(0)\bigr)$ we say that $\phe$ is {\sl tangent} to~$(p;v)$.
 
In [AP] we showed that $\phe$ is a segment of geodesic complex curve iff it is
a holomorphic solution of the system
$$(\phe'')^\alpha+A(\phe')^\alpha=-\Gamma^\alpha_{;i}(\phe;\phe')(\phe')^i,
        \qquad\alpha=1,\ldots,n,\neweq\eqdsei$$
$$G_{\alpha\beta}(\phe;\phe')(\phe''+A\phe')^\beta=\bigl[G_{i;\alpha}(\phe;
        \phe')-G_{\alpha;i}(\phe;\phe')\bigr](\phe')^i,\qquad\alpha=1,\ldots,n,
        \neweq\eqdsette$$
where the prime stands for $\de/\de\zeta$, and $A\colon U\to\C$ is the
function
$$A(\zeta)=-{2\smb\zeta\over1-|\zeta|^2}.$$
 
As we shall see later on, the main amount of informations is contained in
equation~\eqdsei. For the moment, however, let us discuss equation \eqdsette\
a bit.
 
Let $\phe\colon U_r\to M$ be a holomorphic solution of equation \eqdsei. Then,
putting \eqdsei\ into~\eqdsette, we find that $\phe$ is a segment of geodesic
complex curve iff
$$G_{\alpha\beta}(\phe;\phe')\Gamma^\beta_{;i}(\phe;\phe')(\phe')^i=
        \bigl[G_{\alpha;i}(\phe;\phe')-G_{;\alpha}(\phe;\phe')\bigr](\phe')^i,$$
that is iff along the curve $\phe$ we have
$$\bigl[G_{i;\alpha}-G_{\alpha;i}+G_{\alpha\beta}\Gamma^\beta_{;i}\bigr]v^i=0,
        \qquad\alpha=1,\ldots,n.\neweq\eqddieci$$
In a more symmetric way, following [Ru2] we may introduce the {\sl torsion
tensor}
$$T_{\alpha i\smb\mu}=(G_{i\smb\mu;\alpha}-G_{i\smb\mu\beta}\Gamma^\beta
        _{;\alpha})-(G_{\alpha\smb\mu;i}-G_{\alpha\smb\mu\beta}\Gamma^\beta
        _{;i});$$
it is a (3,0)-tensor defined on $(T^{1,0}M)_0$. Then \eqddieci\ is equivalent to
$$T_{\alpha i\smb\mu}\,\bar{v^\mu}v^i=0,\qquad \alpha=1,\ldots,n.$$
 
If $G(p;v)=g_{\alpha\smb\beta}(p)\,v^\alpha\bar{v^\beta}$ is a standard
hermitian metric, then $G_{\alpha\beta}\equiv0$ and \eqddieci\ reduces to
$${\de g_{\alpha\smb\mu}\over\de z^i}={\de g_{i\smb\mu}\over\de z^
        \alpha},$$
that is to the usual K\"ahler condition. For this reason, a strongly
pseudoconvex smooth complex Finsler metric satisfying \eqddieci\ will be said
{\sl K\"ahler.} Summing up, we have proved
 
\newthm Proposition \dT: Let $F$ be a strongly pseudoconvex smooth complex
Finsler metric on a manifold $M$. Then a holomorphic solution $\phe$ of
$\eqdsei$ is a segment of geodesic complex curve iff $F$ is K\"ahler
along~$\phe$.
 
It is possible to write \eqddieci\ in still another way. Set
$\Gamma^\alpha_{\beta;i}=(\Gamma^\alpha_{;i})_\beta$;
then, using \eqde, \eqda\ and \eqdc,
$$\Gamma^\alpha_{\beta;i}=-G^{\gamma\smb\mu}G^{\alpha\smb\nu}
        G_{\gamma\smb\nu\beta}G_{\smb\mu;i}+G^{\alpha\smb\mu}G_{\beta\smb\mu
        ;i},\neweq\eqtu$$
and so
$$\displaylines{G_\alpha\Gamma^\alpha_{\beta;i}=G_{\beta;i}-G_{\beta\gamma}
        \Gamma^\gamma_{;i},\cr
        G_{\alpha\smb\mu}\Gamma^\alpha_{\beta;i}=G_{\beta\smb\mu;i}-G_{\beta
        \smb\mu\gamma}\Gamma^\gamma_{;i},\cr
        G_\alpha\Gamma^\alpha_{i;\beta}v^i=G_{i;\beta}v^i.\cr}$$
Then
$$T_{\alpha i\smb\mu}=G_{\beta\smb\mu}(\Gamma^\beta_{i;\alpha}-\Gamma^\beta
        _{\alpha;i}),\neweq\eqtK$$
$$[G_{i;\alpha}-G_{\alpha;i}+G_{\alpha\beta}\Gamma^\beta_{;i}]v^i=
        G_\beta(\Gamma^\beta_{i;\alpha}-\Gamma^\beta_{\alpha;i})v^i,
        \neweq\eqduuno$$
and \eqddieci\ is equivalent to $G_\beta(\Gamma^\beta_{i;\alpha}-\Gamma^\beta
_{\alpha;i})v^i=0$ for $\alpha=1,\ldots,n$. We remark that when $G$ is a
hermitian metric then
$$\Gamma^\alpha_{\beta;i}=g^{\alpha\smb\mu}\,{\de g_{\beta\smb\mu}\over\de
        z^i},$$
and so they are the coefficients of the Cartan-Chern connection associated to
the hermitian metric. In particular, then, \eqtK\ shows that in this case $T$
actually coincides with the torsion tensor of the connection.
 
But let us now return to equation \eqdsei\ and holomorphic curvature.
We shall say that a holomorphic curve $\phe\colon U_r\to M$ {\sl realizes the
holomorphic curvature at~$0$} if
$$K(\phe^*G)(0)=K_F\bigl(\phe(0);\phe'(0)\bigr).$$
More generally, $\phe$ {\sl realizes the holomorphic curvature at}
$\zeta_0\in U_r$ if $\phe\circ\gamma_{\zeta_0}$ realizes it at~$0$, where
$$\gamma_{\zeta_0}(\zeta)={\zeta+\zeta_0\over1+\bar{\zeta_0}\zeta}$$
is the unique automorphism of $U$ sending the origin to $\zeta_0$ with
positive derivative at~$0$.
 
\newthm Proposition \dtre: Let $F$ be a strongly pseudoconvex smooth complex
Finsler metric on a manifold $M$, and let $\phe\colon U_r\to M$ be a
holomorphic solution of $\eqdsei$. Then
{\smallskip
\itm{(i)}$\phe$ realizes the holomorphic curvature at every point of $U_r$;
\itm{(ii)}if $\bigl(\phe(0);\phe'(0)\bigr)\in S^{1,0}M$, then $\phe$ is a
segment of infinitesimal complex geodesic for~$F$.}
 
\pf (i) By Lemma \duno\ and \eqdtre, a holomorphic $\phe\colon U_r\to M$
realizes the holomorphic curvature at $0$ iff
$$\eta^\alpha=-\lambda^2G^{\alpha\smb\beta}(p;\xi)B_{\smb\beta;i}(p;\xi)\,\xi^i
        +c\,\xi^\alpha,\qquad\alpha=1,\ldots,n,\neweq\eqdotto$$
where $p=\phe(0)$, $v=\phe'(0)=\lambda\xi$ with $\xi\in S^{1,0}_pM$, $\eta=
\phe''(0)$ and $c\in\C$. Since, by \eqda, we have
$$G^{\alpha\smb\beta}B_{\smb\beta;i}=\Gamma^\alpha_{;i}-G_{;i}\xi^\alpha,$$
it follows that \eqdotto\ is equivalent to
$$\eta^\alpha=-\Gamma^\alpha_{;i}(p;v)v^i+c_1v^\alpha,\qquad\alpha=1,\ldots,n,
        \neweq\eqdnove$$
with a possibly different $c_1\in\C$. But \eqdnove\ with $c_1=0$ is just
\eqdsei\ evaluated in~0; so a holomorphic solution of \eqdsei\ realizes the
holomorphic curvature at the origin.
 
Now take $\zeta_0\in U_r$, and set $\psi=\phe\circ\gamma_{\zeta_0}$. Then
$$\displaylines{\psi(0)=\phe(\zeta_0);\cr
        \psi'(0)=(1-|\zeta_0|^2)\phe'(\zeta_0);\cr
        \psi''(0)=(1-|\zeta_0|^2)^2\bigl(\phe''(\zeta_0)+A(\zeta_0)
        \phe'(\zeta_0)\bigr).\cr}$$
So $\phe$ realizes the holomorphic curvature at $\zeta_0$ iff
$$\bigl(\phe''(\zeta_0)+A(\zeta_0)\phe'(\zeta_0)\bigr)^\alpha=-\Gamma^\alpha
        _{;i}\bigl(\phe(\zeta_0);\phe'(\zeta_0)\bigr)\bigl(\phe'(\zeta_0)
        \bigr)^i+c_2\bigl(\phe'(\zeta_0)\bigr)^\alpha,$$
and (i) follows.
\smallskip
(ii) Assume \eqdsei\ holds. Then, recalling \eqda, we get
$$\eqalign{G_\alpha(\phe;\phe')(\phe'')^\alpha+AG(\phe;\phe')&=G_\alpha(\phe;
        \phe')\bigl[(\phe'')^\alpha+A(\phe')^\alpha\bigr]\cr
        &=-G_\alpha(\phe;\phe')\Gamma^\alpha_{;i}(\phe;\phe')(\phe')^i=-G_{;i}
        (\phe;\phe')(\phe')^i.\cr}$$
Therefore
$${\de\over\de\zeta}\bigl[G(\phe;\phe')\bigr]=G_{;i}(\phe;\phe')(\phe')^i+
        G_\alpha(\phe;\phe')(\phe'')^\alpha=-AG(\phe;\phe').$$
Now, along the curve $t\mapsto e^{i\theta}t$ we have
$${\de\over\de\zeta}={1\over 2}e^{-i\theta}{d\over dt};$$
therefore $t\mapsto G\bigl(\phe(e^{i\theta}t);\phe'(e^{i\theta}t)\bigr)$ is a
solution of the Cauchy problem
$$\cases{\displaystyle f'(t)={4t\over1-t^2}f(t),\cr
\noalign{\smallskip}
        f(0)=1.\cr}$$
But $f(t)=(1-t^2)^{-2}$ is a solution of the same problem; therefore
$$G\bigl(\phe(e^{i\theta}t);\phe'(e^{i\theta}t)\bigr)\equiv{1\over(1-t^2)^2},$$
and we are done.\qedn
 
So the main point now is to find when \eqdsei\ has a holomorphic solution.
Assume $\phe\colon U_r\to M$ is such a solution, and apply
$\de/\de\smb\zeta$ to \eqdsei. We get
$${2\over(1-|\zeta|^2)^2}(\phe')^\alpha=\Gamma^\alpha_{;i\smb\j}(\phe;\phe')
        (\phe')^i\bar{(\phe')^j}+\Gamma^\alpha_{\smb\beta;i}(\phe;\phe')
        (\phe')^i\bar{(\phe'')^\beta},\neweq\eqdudue$$
where, as before, $\Gamma^\alpha_{;i\smb\j}=(\Gamma^\alpha_{;i})_{;\smb\j}$
and $\Gamma^\alpha_{\smb\beta;i}=(\Gamma^\alpha_{;i})_{\smb\beta}$.
 
Now, \eqdf\ and \eqdg\ yield
$$\Gamma^\alpha_{\smb\beta;i}\,\bar{v^\beta}=G^{\alpha\smb\mu}_{\smb\beta}
        G_{\smb\mu;i}\,\bar{v^\beta}+G^{\alpha\smb\mu}G_{\smb\beta\smb\mu;i}
        \,\bar{v^\beta}=0.\neweq\eqtnmez$$
So we can replace $\phe''$ by $\phe''+A\phe'$ in \eqdudue\ and, calling in
\eqdsei\ again, we obtain
$${2\over(1-|\zeta|^2)^2}(\phe')^\alpha=\Gamma^\alpha_{;i\smb\j}(\phe;\phe')
        (\phe')^i\bar{(\phe')^j}-\Gamma^\alpha_{\smb\beta;i}(\phe;\phe')
        \Gamma^{\smb\beta}_{;\smb\j}(\phe;\phe')(\phe')^i\bar{(\phe')^j}.
        \neweq\eqdutre$$
Now, if $\bigl(\phe(0);\phe'(0)\bigr)\in S^{1,0}M$, then
Proposition~\dtre.(ii) yields
$$G\bigl(\phe(\zeta);\phe'(\zeta)\bigr)={1\over(1-|\zeta|^2)^2};$$
therefore --- setting $v=\phe'(\zeta)$ --- \eqdutre\ becomes
$$\bigl[\Gamma^\alpha_{;i\smb\j}-\Gamma^\alpha_{\smb\beta;i}\Gamma^{\smb\beta}
        _{;\smb\j}\bigr]v^i\bar{v^j}=2G\,v^\alpha,\qquad\alpha=1,\ldots,n.
        \neweq\eqduquattro$$
So \eqduquattro\ is a necessary condition for \eqdsei\ to have a holomorphic
solution. The interesting fact is that it is sufficient too:
 
\newthm Theorem \dquattro: Let $F$ be a strongly pseudoconvex smooth complex
Finsler metric on a manifold $M$. Then the Cauchy problem
$$\cases{(\phe'')^\alpha+A(\phe')^\alpha=-\Gamma^\alpha_{;i}(\phe;\phe')
        (\phe')^i&for $\alpha=1,\ldots,n$,\cr
        \noalign{\smallskip}
        \phe(0)=p,\quad\phe'(0)=v_0,\cr}\neweq\eqducinque$$
admits a holomorphic solution for all $(p;v_0)\in S^{1,0}M$ iff
$\eqduquattro$ holds. Furthermore, the solution, if exists, is unique.
 
\pf We have already proved one direction; so assume \eqduquattro\ holds.
 
For any $e^{i\theta}\in{\bf S}^1$, let consider the Cauchy problem
$$\cases{\bigl(\ddot g(t)\bigr)^\alpha=-A(t)\bigl(\dot g(t)\bigr)^\alpha
        -\Gamma^\alpha_{;i}\bigl(g(t);\dot g(t)\bigr)\bigl(\dot g(t)\bigr)^i,
        &for $\alpha=1,\ldots,n$,\cr
        \noalign{\smallskip}
        g(0)=p,\quad \dot g(0)=e^{i\theta}v_0.\cr}\neweq\eqdusei$$
The standard ODE theory provides us with an $\eps>0$ and uniquely determined
maps $g_{e^{i\theta}}\colon(-\eps,\eps)\to M$ solving \eqdusei. Define $\phe
\colon U_\eps\to M$ by
$$\phe(\zeta)=g_{\zeta/|\zeta|}(|\zeta|),\neweq\eqdusette$$
and assume for a moment that $\phe$ is holomorphic. Since, writing $\zeta=te^{
i\theta}$, we have
$${\de\over\de\zeta}=-{ie^{-i\theta}\over2t}\left[{\de\over\de\theta}+it{\de
        \over\de t}\right]\qquad\hbox{and}\qquad{\de\over\de\smb\zeta}={ie^{i
        \theta}\over2t}\left[{\de\over\de\theta}-it{\de\over\de t}\right],
        \neweq\eqduotto$$
it follows that
$${\de\phe\over\de\zeta}(\zeta)={\smb\zeta\over|\zeta|}\dot g_{\zeta/|\zeta|}
        (|\zeta|),$$
and thus $\phe$ is a holomorphic solution of \eqducinque. In conclusion, we
must prove that, assuming \eqduquattro, the map $\phe$ defined by \eqdusette\
is holomorphic. Note that, since a holomorphic map is uniquely determined by
its restriction to the real axis, the uniqueness statement for \eqdusei\
implies that $\phe$~is the unique possible holomorphic solution of \eqducinque.
 
First of all, set $f_0(t)=\hbox{tanh}\,t$ and
$$\sigma_\theta(t)=g_{e^{i\theta}}(\hbox{tanh}\,t).\neweq\eqduottomez$$
Then
$$\dot\sigma_\theta=(f'_0)(\dot g_{e^{i\theta}}\circ f_0)\qquad\hbox{and}
        \qquad\ddot\sigma_\theta=(f'_0)^2[(\ddot g_{e^{i\theta}}+A\dot
        g_{e^{i\theta}})\circ f_0];$$
so $\sigma_\theta$ satisfies
$$\cases{\ddot\sigma_\theta^\alpha=-\Gamma^\alpha_{;i}(\sigma_\theta;\dot
        \sigma_\theta)\,\dot\sigma_\theta^i,&for $\alpha=1,\ldots,n$,\cr
        \sigma_\theta(0)=p,\quad\dot\sigma_\theta(0)=e^{i\theta}v_0.\cr}$$
Set $h=G(\sigma_\theta;\dot\sigma_\theta)$. Then
$$\eqalign{h'&=G_{;i}(\sigma_\theta;\dot\sigma_\theta)\,\dot\sigma_\theta^i
        +G_\alpha(\sigma_\theta;\dot\sigma_\theta)\,\ddot\sigma_\theta^\alpha\cr
        &=G_{;i}(\sigma_\theta;\dot\sigma_\theta)\,\dot\sigma_\theta^i
        +G_{\alpha\smb\beta}(\sigma_\theta;\dot\sigma_\theta)\,\ddot\sigma_
        \theta^\alpha\bar{\dot\sigma_\theta^\beta}\cr
        &=G_{;i}(\sigma_\theta;\dot\sigma_\theta)\,\dot\sigma_\theta^i
        -G_{\smb\mu;i}(\sigma_\theta;\dot\sigma_\theta)\,\dot\sigma_\theta^i
        \bar{\dot\sigma_\theta^\mu}=0.\cr}$$
So $h(t)\equiv h(0)=1$, and the curve $\sigma_\theta$ lifts to a curve
$\tilde\sigma_\theta=d\sigma_\theta=(\sigma_\theta;\dot\sigma_\theta)$ in $S^{
1,0}M$.
 
Now, we define a global vector field $X\in\Gamma\Bigl(T^{1,0}(S^{1,0}M)\Bigr)$
by setting
$$X_{\tilde v}=v^i{\de\over\de z^i}-\Gamma^\alpha_{;i}\,v^i{\de\over\de
        v^\alpha},\neweq\eqdunove$$
where $(z^1,\ldots,z^n;v^1,\ldots,v^n)$ are the local coordinates of $\tilde
v\in(T^{1,0}M)_0$. It is not difficult to check that $X$ is globally defined,
and that for $\tilde v\in S^{1,0}
M$ it is actually true that $X_{\tilde v}\in T^{1,0}_{\tilde v}(S^{1,0}M)$,
as claimed.
 
To proceed, we need to recall a basic fact of complex differential geometry.
Let $N$ be a complex manifold, of complex dimension~$m$. If we consider $N$
with its real structure, then $TN$ is a $(4m)$-dimensional real vector bundle
on $N$ endowed with a complex structure~$J$. If we denote by $T^cN$ its
complexification, then $T^{1,0}N$ is the $i$-eigenspace of $J$, and the
canonical isomorphism $T^{1,0}N\to TN$ is given by
$$Y\mapsto Y^o=Y+\bar Y,$$
where $\bar Y$ is the complex conjugate of $Y$ in $T^cN$. In particular,
then,
$$JY^o=i(Y-\bar Y).$$
 
The aim of this observation is that, by construction, $\tilde\sigma_\theta$ is
the integral curve in $S^{1,0}M$ of the vector field $X^o$ starting at
$(p;e^{i\theta}v_0)$. If we denote by $e^{tX^o}$ the local one-parameter group
of diffeomorphisms induced by $X^o$ on $S^{1,0}M$, we may then write
$$\tilde\sigma_\theta(t)=e^{tX^o}(e^{i\theta}\tilde v_0)\qquad\hbox{and}
        \qquad\sigma_\theta(t)=\pi(e^{tX^o}e^{i\theta}\tilde v_0),
        \neweq\eqdunovemez$$
where $\pi\colon S^{1,0}M\to M$ is the canonical projection, and $\tilde v_0=
(p;v_0)\in S^{1,0}M$.
 
We need another vector field on $S^{1,0}M$. The map $\bigl(e^{i\theta},(p;v)
\bigr)\mapsto(p;e^{i\theta}v)$ is a one-parameter group of diffeomorphisms of $S
^{1,0}M$; therefore it is induced by a vector field $Z$, namely
$$Z=iv^\alpha{\de\over\de v^\alpha}\in\Gamma\Bigl(T^{1,0}(S^{1,0}M)
        \Bigr);$$
note that $\pi_*(Z)=0$. Then \eqdunovemez\ becomes
$$\tilde\sigma_\theta(t)=e^{tX^o}e^{\theta Z^o}\tilde v_0\qquad\hbox{and}
        \qquad\sigma_\theta(t)=\pi(e^{tX^o}e^{\theta Z^o}\tilde v_0).
        \neweq\eqdventi$$
 
Now, we need to compute
$$[X^o,JX^o]=i[X+\bar X,X-\bar X]=-2i[X,\bar X],$$
and
$$[X^o,Z^o]=[X+\bar X,Z+\bar Z]=[X,Z]^o+[X,\bar Z]^o.$$
Using local coordinates we find
$$\eqalign{[X,Z]&=-i\Gamma^\alpha_{;j}\,v^j{\de\over\de v^\alpha}-iv^\beta
        \left[{\de\over\de z^\beta}-\Gamma^\alpha_{\beta;j}\,v^j{\de\over\de
        v^\alpha}-\Gamma^\alpha_{;\beta}{\de\over\de v^\alpha}\right]\cr
        &=-i\left[v^\beta{\de\over\de z^\beta}-\Gamma^\alpha_{;j}\,v^j
        {\de\over\de v^\alpha}\right]=-iX,\cr}\neweq\eqbruno$$
because, by \eqdf\ and \eqdg, $\Gamma^\alpha_{\beta;j}\,v^\beta=\Gamma^\alpha_
{;j}$. It is clear by the definitions and \eqtnmez\ that $[X,\bar Z]=0$;
finally,
$$\eqalign{[X,\bar X]&=-\left[\Gamma^{\smb\alpha}_{;\smb hj}-\Gamma^\beta_{;j}
        \Gamma^{\smb\alpha}_{\beta;\smb h}\right]v^j\bar{v^h}{\de\over\de
        \bar{v^\alpha}}+\left[\Gamma^\alpha_{;h\smb\j}-\Gamma^{\smb\beta}_
        {;\smb\j}\Gamma^\alpha_{\smb\beta;h}\right]\bar{v^j}v^h{\de\over\de
        v^\alpha}\cr
        &=2\left(v^\alpha{\de\over\de v^\alpha}-\bar{v^\alpha}{\de\over\de
        \bar{v^\alpha}}\right)=-2iZ^o,\cr}$$
where we used \eqduquattro\ on $S^{1,0}M$. So we get
$$[X^o,JX^o]=-4Z^o\qquad\hbox{and}\qquad[X^o,Z^o]=-JX^o.$$
 
Now fix $\tau>0$, and set $\tilde v_\tau=e^{\tau X^o}\tilde v_0$. Put
$$u(t)=e^{tX^o}_*Z^o_{e^{-tX^{\!o}}\tilde v_\tau}\in T_{\tilde v_\tau}(
        S^{1,0}M).$$
Then
$${du\over dt}(t)={d\over dt}\bigl(e^{tX^o}_*Z^o_{e^{-tX^{\!o}}\tilde v_\tau}
        \bigr)=-e^{tX^o}_*\{{\cal L}_{X^o}Z^o\}_{e^{-tX^{\!o}}\tilde v_\tau}
        =e^{tX^o}_*(JX^o)_{e^{-tX^{\!o}}\tilde v_\tau},$$
where ${\cal L}_{X^o}$ is the Lie derivative, and
$${d^2u\over dt^2}(t)={d\over dt}\left(e^{tX^o}_*(JX^o)_{e^{-tX^{\!o}}\tilde
        v_\tau}\right)=-e^{tX^o}_*\{{\cal L}_{X^o}(JX^o)\}_{e^{-tX^{\!o}}
        \tilde v_\tau}=4e^{tX^o}_*Z^o_{e^{-tX^{\!o}}\tilde v_\tau}.$$
In other words, $u(t)$ is a solution of the Cauchy problem
$$\cases{\ddot u=4u,\cr
        \noalign{\smallskip}
        u(0)=Z^o_{\tilde v_\tau},\quad\dot u(0)=(JX^o)_{\tilde v_\tau}.\cr}$$
Therefore
$$u(t)={1\over4}e^{2t}\bigl(2Z^o_{\tilde v_\tau}+(JX^o)_{\tilde v_\tau}\bigr)
        +{1\over 4}e^{-2t}\bigl(2Z^o_{\tilde v_\tau}-(JX^o)_{\tilde v_\tau}
        \bigr),$$
and, in particular,
$$\pi_*e^{\tau X^o}_*Z^o_{\tilde v_0}=\pi_*u(\tau)={e^{2\tau}-e^{-2\tau}\over4}
        \pi_*(JX^o)_{\tilde v_\tau}={e^{2\tau}-e^{-2\tau}\over4}J\pi_*X^o_
        {\tilde v_\tau}.\neweq\eqdduno$$
 
We are almost done. Recalling \eqdusette, \eqduotto, \eqduottomez\ and
\eqdventi, it is clear that we should prove that
$$\left.{\de\over\de\theta}\pi\left(e^{(\atanh t)X^o}e^{\theta
        Z^o}\tilde v_0\right)\right|_{\theta=0}=tJ\left.{\de\over\de t}\pi
        \left(e^{(\atanh t)X^o}e^{\theta Z^o}\tilde v_0\right)
        \right|_{\theta=0},$$
where we may take $\theta=0$ because $\tilde v_0$ is generic. Let us
compute; using \eqdduno\ we get
$$\left.{\de\over\de\theta}\pi\left(e^{(\atanh t)X^o}e^{\theta
        Z^o}\tilde v_0\right)\right|_{\theta=0}=\pi_*e^{(\atanh t)X^o}_*
        Z^o_{\tilde v_0}={t\over1-t^2}J\pi_*X^o_{\tilde v_{\atanh t}},$$
whereas
$$\left.{\de\over\de t}\pi\left(e^{(\atanh t)X^o}e^{\theta
        Z^o}\tilde v_0\right)\right|_{\theta=0}={\de\over\de t}\pi\left(
        e^{(\atanh t)X^o}\tilde v_0\right)={1\over1-t^2}
        \pi_*X^o_{\tilde v_{\atanh t}},$$
and the proof is complete.\qedn
 
So we have found a necessary and sufficient condition for the existence of
segments of geodesic complex curves:
 
\newthm Corollary \dsei: Let $F$ be a strongly pseudoconvex smooth complex
Finsler metric on a manifold $M$. Then:
{\smallskip
\itm{(i)}if $\eqduquattro$ holds, then for any $(p;\xi)\in S^{1,0}M$ there is
a segment of infinitesimal complex geodesic tangent to $(p;\xi)$;
\itm{(ii)}there exists a (unique) segment of geodesic complex curve tangent to
$(p;\xi)$ for any $(p;\xi)\in S^{1,0}M$ iff $F$ is K\"ahler and $\eqduquattro$
holds.}
 
\pf (i) Theorem \dquattro\ and Proposition \dtre.(ii).
 
(ii) In [AP] it is shown that a segment of geodesic complex curve is a
holomorphic solution of the system \eqdsei--\eqdsette. The assertion then
follows from Theorem \dquattro\ and Proposition~\dT.\qedn
 
A natural question now is whether the completeness of the metric $F$ ---
together with K\"ahler and \eqduquattro\ --- would imply the existence of
geodesic complex curves defined on the whole unit disk~$U$. The answer is
positive, but for the proof we beforehand need a discussion of the geometrical
meaning of \eqduquattro.
 
Thanks to Proposition~\dtre.(i), we know that a holomorphic solution of
\eqdsei\ realizes the holomorphic curvature and it is an isometry from the
Poincar\'e metric to $F$; in particular, thus, the holomorphic curvature
along the curve should be~$-4$. This suggests to look for a connection between
\eqduquattro\ and the holomorphic curvature; and indeed the next result shows
that the connection is provided by a sort of simmetry condition on the
curvature. Analogously to the tensor $T_{\alpha i\smb\mu}$ previously
introduced, set
$$\eqalign{H_{\alpha i\smb\mu\smb\j}&=G_{\tau\smb\mu}(\Gamma^\tau_{i;\alpha}
        -\Gamma^\tau_{\alpha;i})_{;\smb\j}+G_{\tau\smb\mu i}\Gamma^\tau
        _{;\alpha\smb\j}-G_{\tau\smb\mu\alpha}\Gamma^\tau_{;i\smb\j}\cr
        &=(G_{\tau\smb\mu}\Gamma^\tau_{;\alpha\smb\j})_i-(G_{\tau\smb\mu}
        \Gamma^\tau_{;i\smb\j})_\alpha;\cr}$$
it is a (4,0)-tensor on $(T^{1,0}M)_0$. Note that
$$H_{\alpha i\smb\mu\smb\j}\,\bar{v^\mu}v^i\bar{v^j}=\bigl[G_\tau(\Gamma^
        \tau_{i;\alpha}-\Gamma^\tau_{\alpha;i})_{;\smb\j}-G_{\tau\alpha}
        \Gamma^\tau_{;i\smb\j}\bigr]v^i\bar{v^j}.\neweq\eqdH$$
Then
 
\newthm Theorem \dC: Let $F$ be a strongly pseudoconvex smooth complex Finsler
metric on a manifold $M$. Then $\eqduquattro$ holds iff $K_F\equiv-4$ and
$$H_{\alpha i\smb\mu\smb\j}\,\bar{v^\mu}v^i\bar{v^j}=0,\qquad\alpha=1,
        \ldots,n.\neweq\eqdcro$$
 
\pf We start by showing that \eqduquattro\ implies $K_F\equiv-4$.
Indeed, take $(p;\xi)\in S^{1,0}M$. Then, recalling \eqdf, \eqda, and \eqdh,
we get
$$G_\alpha\Gamma^\alpha_{\smb\beta;i}=G_\alpha G^{\alpha\smb\mu}_{\smb\beta}
        G_{\smb\mu;i}+G_\alpha G^{\alpha\smb\mu}G_{\smb\mu\smb\beta;i}=0.$$
Therefore \eqduquattro\ yields
$$K_F(p;\xi)=-2G_\alpha\Gamma^{\alpha}_{;i\smb\j}\,\xi^i\bar{\xi^j}=-4G_\alpha
        \xi^\alpha-G_\alpha\Gamma^\alpha_{\smb\beta;i}\Gamma^{\smb\beta}
        _{\smb\j}\,\xi^i\bar{\xi^j}=-4G_\alpha\xi^\alpha=-4.$$
 
From now on we shall assume $K_F\equiv-4$; we ought to prove that in this case
\eqduquattro\ is equivalent to \eqdcro, that is, by \eqdH, to
$$G_\tau(\Gamma^\tau_{i;\alpha}-\Gamma^\tau_{\alpha;i})_{;\smb\j}\,v^i
        \bar{v^j}=G_{\tau\alpha}\Gamma^\tau_{;i\smb\j}\,v^i\bar{v^j},\qquad
        \alpha=1,\ldots,n.\neweq\eqdcroce$$
 
By \eqdcinquemez, $K_F\equiv-4$ is equivalent to
$$G_\beta\Gamma^\beta_{;i\smb\j}\,v^i\bar{v^j}=2G^2.$$
Differentiating with respect to $\bar{v^\nu}$ we get
$$4GG_{\smb\nu}=\bigl[G_{\beta\smb\nu}\Gamma^\beta_{;i\smb\j}+G_\beta
        \Gamma^\beta_{\smb\nu;i\smb\j}\bigr]v^i\bar{v^j}+G_\beta\Gamma^\beta
        _{;i\smb\nu}\,v^i;$$
multiplying by $G^{\alpha\smb\nu}$, and recalling \eqda, we obtain
$$4G\,v^\alpha=\bigl[\Gamma^\alpha_{;i\smb\j}+G^{\alpha\smb\nu}G_\beta
        \Gamma^\beta_{\smb\nu;i\smb\j}\bigr]v^i\bar{v^j}+G^{\alpha\smb\nu}
        G_\beta\Gamma^\beta_{;i\smb\nu}v^i.\neweq\eqdmenouno$$
Now,
$$\Gamma^\beta_{;i\smb\j}=G^{\beta\smb\mu}G_{\smb\mu;i\smb\j}-G^{\beta\smb\tau}
        G_{\sigma\smb\tau;\smb\j}\Gamma^\sigma_{;i};\neweq\eqdgauno$$
$$\Gamma^\beta_{\smb\nu;i\smb\j}=G^{\beta\smb\mu}_{\smb\nu}G_{\smb\mu;
        i\smb\j}+G^{\beta\smb\mu}G_{\smb\mu\smb\nu;i\smb\j}-G^{\beta\smb\tau}
        _{\smb\nu}G_{\sigma\smb\tau;\smb\j}\Gamma^\sigma_{;i}-G^{\beta\smb\tau}
        G_{\sigma\smb\tau\smb\nu;\smb\j}\Gamma^\sigma_{;i}-G^{\beta\smb\tau}
        G_{\sigma\smb\tau;\smb\j}\Gamma^\sigma_{\smb\nu;i}.$$
Therefore, using \eqdf, \eqde, \eqdb, \eqda, \eqdc\ and \eqdh, we get
$$\eqalign{G^{\alpha\smb\nu}G_\beta\Gamma^\beta_{\smb\nu;i\smb\j}&=-G^{\alpha
        \smb\nu}G_{\sigma;\smb\j}\Gamma^\sigma_{\smb\nu;i}=-G^{\alpha\smb\nu}
        G_{\sigma;\smb\j}(G^{\sigma\smb\mu}_{\smb\nu}G_{\smb\mu;i}+G^{\sigma
        \smb\mu}G_{\smb\mu\smb\nu;i})\cr
        &=-\Gamma^{\smb\mu}_{;\smb\j}G^{\alpha\smb\nu}G_{\smb\mu\smb\nu;i}
        +\Gamma^{\smb\gamma}_{;\smb\j}G^{\alpha\smb\nu}G^{\delta\smb\mu}
        G_{\delta\smb\gamma\smb\nu}G_{\smb\mu;i}\cr
        &=-\Gamma^{\smb\mu}_{;\smb\j}(G^{\alpha\smb\nu}G_{\smb\mu\smb\nu;i}
        +G^{\alpha\smb\nu}_{\smb\mu}G_{\smb\nu;i})=-\Gamma^\alpha_{\smb\mu;i}
        \Gamma^{\smb\mu}_{\smb\j}.\cr}$$
So \eqdmenouno\ becomes
$$4G\,v^\alpha=\bigl[\Gamma^\alpha_{;i\smb\j}-\Gamma^\alpha_{\smb\mu;i}
        \Gamma^{\smb\mu}_{;\smb\j}\bigr]v^i\bar{v^j}+G^{\alpha\smb\nu}G_\beta
        \Gamma^\beta_{;i\smb\nu}\,v^i.\neweq\eqdmenounob$$
Now, in \eqduseibis\ we showed that
$$G^{\alpha\smb\nu}G_\beta\Gamma^\beta_{;i\smb\nu}\,v^i=G^{\alpha\smb\nu}\bigl[
        G_{;i\smb\nu}-\Gamma^{\smb\tau}_{;\smb\j}G_{\smb\tau;i}\bigr]v^i.$$
Since, by \eqdb\ and \eqdg,
$$\Gamma^{\smb\tau}_{\smb\j;\smb\nu}\,\bar{v^j}=G^{\sigma\smb\tau}G_{\sigma
        \smb\j;\smb\nu}\,\bar{v^j}=\Gamma^{\smb\tau}_{;\smb\nu},$$
we get
$$\eqalign{G^{\alpha\smb\nu}G_\beta\Gamma^\beta_{;i\smb\nu}\,v^i&=G^{\alpha
        \smb\nu}\bigl[G_{\smb\j;i\smb\nu}-G_{\smb\tau;i}\Gamma^{\smb\tau}
        _{\smb\j;\smb\nu}\bigr]v^i\bar{v^j}\cr
        &=G^{\alpha\smb\nu}\bigl[G_{\smb\nu;i\smb\j}-G_{\smb\tau;i}\Gamma^
        {\smb\tau}_{\smb\nu;\smb\j}\bigr]v^i\bar{v^j}+G^{\alpha\smb\nu}
        \bigr[(G_{\smb\j;\smb\nu}-G_{\smb\nu;\smb\j})_{;i}-G_{\smb\tau;i}
        (\Gamma^{\smb\tau}_{\smb\j;\smb\nu}-\Gamma^{\smb\tau}_{\smb\nu;\smb\j})
        \bigr]v^i\bar{v^j}.\cr}$$
Now \eqdc\ yields
$$\eqalign{G^{\alpha\smb\nu}G_{\smb\tau;i}\Gamma^{\smb\tau}_{\smb\nu;\smb\j}&=
        G^{\alpha\smb\nu}G_{\sigma\smb\nu;\smb\j}\Gamma^\sigma_{;i}-G^{\alpha
        \smb\nu}G_{\delta\smb\mu\smb\nu}\Gamma^{\smb\mu}_{;\smb\j}
        \Gamma^\delta_{;i}\cr
        &=G^{\alpha\smb\nu}G_{\sigma\smb\nu;\smb\j}\Gamma^\sigma_{;i}+
        \Gamma^\alpha_{\smb\mu;i}\Gamma^{\smb\mu}_{;\smb\j}-G^{\alpha\smb\nu}
        G_{\smb\mu\smb\nu;i}\Gamma^{\smb\mu}_{;\smb\j};\cr}$$
so, by \eqdgauno,
$$G^{\alpha\smb\nu}(G_{\smb\nu;i\smb\j}-G_{\smb\tau;i}\Gamma^{\smb\tau}
        _{\smb\nu;\smb\j})=(\Gamma^\alpha_{;i\smb\j}
        -\Gamma^\alpha_{\smb\mu;i}\Gamma^{\smb\mu}_{;\smb\j})
        +G^{\alpha\smb\nu}G_{\smb\mu\smb\nu;i}\Gamma^{\smb\mu}_{;\smb\j}.$$
Summing up, we have found
$$\eqalign{G^{\alpha\smb\nu}G_\beta\Gamma^\beta_{;i\smb\nu}v^i=&\bigl[\Gamma^
        \alpha_{;i\smb\j}-\Gamma^\alpha_{\smb\mu;i}\Gamma^{\smb\mu}_{;\smb\j}
        \bigr]v^i\bar{v^j}\cr
        &\quad+G^{\alpha\smb\nu}\bigl[G_{\smb\tau\smb\nu;i}\Gamma^{\smb\tau}
        _{;\smb\j}+(G_{\smb\j;\smb\nu}-G_{\smb\nu;\smb\j})_{;i}-G_{\smb\tau;i}
        (\Gamma^{\smb\tau}_{\smb\j;\smb\nu}-\Gamma^{\smb\tau}_{\smb\nu;\smb\j})
        \bigr]v^i\bar{v^j}.\cr}$$
For the moment, set
$$T_{\smb\nu}=[G_{\smb\j;\smb\nu}-G_{\smb\nu;\smb\j}+G_{\smb\tau\smb\nu}\Gamma
        ^{\smb\tau}_{;\smb\j}\bigr]\bar{v^j}=T_{\smb\nu\smb\j\mu}\,v^\mu
        \bar{v^j};$$
recall that \eqduuno\ says that
$$T_{\smb\nu}=G_{\smb\tau}(\Gamma^{\smb\tau}_{\smb\j;\smb\nu}-\Gamma^{\smb\tau}
        _{\smb\nu;\smb\j})\bar{v^j}.$$
Therefore
$$\displaylines{(G_{\smb\j;\smb\nu}-G_{\smb\nu;\smb\j})_{;i}\,\bar{v^j}=(T_
        {\smb\nu})_{;i}-(G_{\smb\tau\smb\nu}\Gamma^{\smb\tau}_{;\smb\j})_{;i}\,
        \bar{v^j}=(T_{\smb\nu})_{;i}-\bigl[G_{\smb\tau\smb\nu;i}\Gamma^{\smb
        \tau}_{;\smb\j}+G_{\smb\tau\smb\nu}\Gamma^{\smb\tau}_{;\smb\j i}\bigr]
        \bar{v^j};\cr
        G_{\smb\tau;i}(\Gamma^{\smb\tau}_{\smb\j;\smb\nu}-\Gamma^{\smb\tau}
        _{\smb\nu;\smb\j})\bar{v^j}=(T_{\smb\nu})_{;i}-G_{\smb\tau}
        (\Gamma^{\smb\tau}_{\smb\j;\smb\nu}-\Gamma^{\smb\tau}
        _{\smb\nu;\smb\j})_{;i}\,\bar{v^j}.\cr}$$
In conclusion, we have shown that
$$G^{\alpha\smb\nu}G_\beta\Gamma^\beta_{;i\smb\nu}v^i=\bigl[\Gamma^\alpha_{;i
        \smb\j}-\Gamma^\alpha_{\smb\mu;i}\Gamma^{\smb\mu}_{;\smb\j}\bigr]v^i
        \bar{v^j}+G^{\alpha\smb\nu}\bigl[G_{\smb\tau}(\Gamma^{\smb\tau}_{\smb\j;
        \smb\nu}-\Gamma^{\smb\tau}_{\smb\nu;\smb\j})_{;i}-G_{\smb\nu\smb\tau}
        \Gamma^{\smb\tau}_{;\smb\j i}\bigr]v^i\bar{v^j}.$$
Recalling \eqdmenounob, we have obtained
$$4G\,v^\alpha=2\bigl[\Gamma^\alpha_{;i\smb\j}-\Gamma^\alpha_{\smb\mu;i}\Gamma
        ^{\smb\mu}_{;\smb\j}\bigr]v^i\bar{v^j}+G^{\alpha\smb\nu}\bigl[G_
        {\smb\tau}(\Gamma^{\smb\tau}_{\smb\j;\smb\nu}-\Gamma^{\smb\tau}_
        {\smb\nu;\smb\j})_{;i}-G_{\smb\nu\smb\tau}\Gamma^{\smb\tau}_{;\smb\j i}
        \bigr]v^i\bar{v^j},$$
and the assertion follows.\qedn
 
If $G(p;v)=g_{\mu\smb\nu}(p)\,v^\mu\bar{v^\nu}$ is an hermitian metric, then
the tensor $H_{\alpha i\smb\mu\smb\j}$ becomes
$$H_{\alpha i\smb\mu\smb\j}=g_{\tau\smb\mu}\,{\de\over\de\bar{z^j}}T^\tau_
        {i\alpha}=R_{i\alpha\smb\mu\smb\j}-R_{\alpha i\smb\mu\smb\j},$$
where $T^\tau_{i\alpha}$ is the torsion of the Chern connection
associated to the hermitian metric, and $R_{\alpha i\smb\mu\smb\j}$ is the
Riemannian curvature tensor of the connection. So \eqdcro\ is equivalent to
$$R_{i\alpha\smb\mu\smb\j}\,v^i\bar{v^\mu}\bar{v^j}=R_{\alpha i\smb\mu\smb\j}\,
        v^i\bar{v^\mu}\bar{v^j}\neweq\eqsim$$
for all $v\in T^{1,0}M$. So \eqdcro\ may be interpreted as a simmetry
condition on a curvature tensor; more precisely, as anticipated in the
introduction, a simmetry condition on the Chern connection induced by the
Finsler metric on the vertical subbundle of $T^{1,0}\bigl(T^{1,0}M\bigr)$.
Finally, we also remark that --- at least in the hermitian case --- \eqdcro\
in particular holds when $\bar\de T\equiv0$.
 
Now that we have an idea of the geometrical meaning of \eqduquattro, we may
return to the study of geodesic complex curves. As anticipated, we are now able
to prove the following
 
\newthm Theorem \dotto: Let $F$ be a strongly pseudoconvex smooth complete
complex Finsler metric on a manifold~$M$. Assume that the holomorphic
curvature of $F$ is identically~$-4$ and that $\eqdcro$ holds. Then for every
$(p;\xi)\in S^{1,0} M$ there is a unique holomorphic solution $\phe\colon U\to
M$ of $\eqdsei$ defined on the whole unit disk~$U$ such that $\phe(0)=p$ and
$\phe'(0)=\xi$.
 
\pf First of all, we remark that the distribution $\C X^o\oplus\C Z^o\subset
T(S^{1,0}M)$ is involutive. Indeed, Theorem~\dC\ shows that
\eqduquattro\ holds, and in the proof of Theorem~\dquattro\ we have already
computed $[X^o,JX^o]=-4Z^o$ and $[X^o,Z^o]=-JX^o$. For the remaining brackets,
using \eqbruno\ we get
$$[X^o,JZ^o]=X^o,\quad[JX^o,Z^o]=X^o,\quad[JX^o,JZ^o]=JX^o,\quad[Z^o,JZ^o]=0.$$
Let $\widetilde N$ denote the integral leaf of this distribution passing
through $(p;\xi)$. From the proof of Theorem~\dquattro\ it follows that
$N=\pi(\widetilde N)\subset M$ is a Riemann surface locally parametrized by
the holomorphic solutions of \eqdsei. In particular, $F$ restricted to $N$ is
a complete hermitian metric of constant Gaussian curvature~$-4$, because of
Proposition~\dtre. Thus there is a unique holomorphic covering map $\psi\colon
U\to N$ which is an isometry between the Poincar\'e metric on~$U$ and $F$
restricted to~$N$ and such that $\psi(0)=p$ and $\psi'(0)=\xi$. But if
$\phe\colon U_\eps\to N$ is the holomorphic solution of \eqdsei\ with $\phe(0)
=p$ and $\phe'(0)=\xi$, then $\phe$ too is an isometry between the Poincar\'e
metric restricted to $U_\eps$ and $F$ restricted to $N$; it follows that $\phe
=\psi|_{U_\eps}$, and $\psi$ is the extension of $\phe$ to the whole~$U$ we
were looking for.\qedn
 
\newthm Corollary \dottobis: Let $F$ be a strongly pseudoconvex smooth
complete complex Finsler metric on a manifold~$M$. Assume that the holomorphic
curvature of $F$ is identically~$-4$ and that $\eqdcro$ holds. Then:
{\smallskip
\itm{(i)} for any $(p;\xi)\in S^{1,0}M$ there is an infinitesimal complex
geodesic tangent to $(p;\xi)$;
\itm{(ii)}if moreover $F$ is K\"ahler, then for any $(p;\xi)\in S^{1,0}M$ there
is a unique geodesic complex curve tangent to~$(p;\xi)$.}
 
\pf Theorems \dC, \dotto\ and Propositions~\dT\ and \dtre.\qedn
 
Actually, we can even get a sort of punctual version of the latter result. As
usual, we need a computation, which by the way clarifies the relationship among
\eqduquattro\ and the two torsion tensors we introduced, $T$ and $H$. For the
sake of simplicity, set
$$\Sigma^\alpha_{;i\smb\j}=\Gamma^\alpha_{;i\smb\j}-\Gamma^\alpha_{\smb\mu;i}
        \Gamma^{\smb\mu}_{;\smb\j}.$$
In particular, \eqduquattro\ becomes $\Sigma^\alpha_{;i\smb\j}\,v^i\bar{v^j}=
2G\,v^\alpha$.
 
\newthm Proposition \dXkahler: Let $F$ be a strongly pseudoconvex smooth
complex Finsler metric on a manifold~$M$. Then
$$H_{\alpha i\smb\mu\smb\j}\,\bar{v^\mu}\,\bar{v^j}=\bar{X}\bigl(T_{\alpha
        i\smb\mu}\,\bar{v^\mu}\bigr)+\bigl(G_{i\sigma}\Sigma^\sigma
        _{;\alpha\smb\j}-G_{\alpha\sigma}\Sigma^\sigma_{;i\smb\j}\bigr)
        \bar{v^j},$$
for all $i$,~$\alpha=1,\ldots,n$, where $\bar{X}$ is the complex conjugate of
the vector field defined in~$\eqdunove$. In particular,
$$H_{\alpha i\smb\mu\smb\j}\,\bar{v^\mu}v^i\bar{v^j}=\bar{X}(T_{\alpha i
        \smb\mu}\,\bar{v^\mu}v^i)-G_{\alpha\sigma}\Sigma^\sigma_{;i\smb\j}\,
        v^i\bar{v^j},\qquad\alpha=1,\ldots,n.$$
 
\pf By \eqtK\ and \eqdunove\
$$\displaylines{T_{\alpha i\smb\mu}\,\bar{v^\mu}=G_\beta(\Gamma^\beta_{i;
        \alpha}-\Gamma^\beta_{\alpha;i}),\cr
\noalign{\smallskip}
        \bar{X}=\bar{v^j}{\de\over\de\bar{z^j}}-\Gamma^{\smb\gamma}_{;\smb\j}
        \,\bar{v^j}{\de\over\de\bar{v^\gamma}}.\cr}$$
Thus
$$\eqalign{\bar{X}(T_{\alpha i\smb\mu}\,\bar{v^\mu})=&G_\beta\bigl[(\Gamma
        ^\beta_{i;\alpha}-\Gamma^\beta_{\alpha;i})_{;\smb\j}-\Gamma^{\smb\gamma}
        _{;\smb\j}(\Gamma^\beta_{i;\alpha}-\Gamma^\beta_{\alpha;i})
        _{\smb\gamma}\bigr]\bar{v^j}\cr
        &\quad+[G_{\beta;\smb\j}-G_{\beta\smb\gamma}\Gamma^{\smb\gamma}
        _{;\smb\j}](\Gamma^\beta_{i;\alpha}-\Gamma^\beta_{\alpha;i})\bar{v^j}.
        \cr}\neweq\eqtdue$$
Now, $G_{\beta\smb\gamma}\Gamma^{\smb\gamma}_{;\smb\j}=G_{\beta;\smb\j}$, and
so the second addendum in \eqtdue\ vanishes. Next, recalling \eqtu\ and the
usual formulas,
$$G_\beta(\Gamma^\beta_{i;\alpha}-\Gamma^\beta_{\alpha;i})_{\smb\gamma}
        =-[G_{i\sigma}\Gamma^\sigma_{\smb\gamma;\alpha}-G_{\alpha\sigma}
        \Gamma^\sigma_{\smb\gamma;i}].$$
Therefore
$$\eqalign{\bar{X}(T_{\alpha i\smb\mu}\,\bar{v^\mu})&=G_\beta(\Gamma^\beta
        _{i;\alpha}-\Gamma^\beta_{\alpha;i})_{;\smb\j}\,\bar{v^j}+
        [G_{i\sigma}\Gamma^\sigma_{\smb\gamma;\alpha}\Gamma^{\smb\gamma}
        _{;\smb\j}-G_{\alpha\sigma}\Gamma^\sigma_{\smb\gamma;i}\Gamma^
        {\smb\gamma}_{;\smb\j}]\bar{v^j}\cr
        &=G_\beta(\Gamma^\beta_{i;\alpha}-\Gamma^\beta_{\alpha;i})_{;\smb\j}
        \,\bar{v^j}-[G_{i\sigma}\Sigma^\sigma_{;\alpha\smb\j}-G_{\alpha
        \sigma}\Sigma^\sigma_{;i\smb\j}]\bar{v^j}+(G_{i\sigma}\Gamma^\sigma
        _{;\alpha\smb\j}-G_{\alpha\sigma}\Gamma^\sigma_{;i\smb\j})\bar{v^j}\cr
        &=H_{\alpha i\smb\mu\smb\j}\,\bar{v^\mu}\,\bar{v^j}-[G_{i\sigma}
        \Sigma^\sigma_{;\alpha\smb\j}-G_{\alpha\sigma}\Sigma^\sigma_{;i\smb\j}
        ]\bar{v^j},\cr}$$
and the assertion follows.\qedn
 
As a corollary we have
 
\newthm Corollary \dXKahlerbis: Let $F$ be a strongly pseudoconvex smooth
complex Finsler metric on a manifold~$M$. Assume that the holomorphic
curvature of $F$ is identically~$-4$ and that $\eqdcro$ holds. Take
$(p_0;\xi_0)\in S^{1,0}M$. Then there is a segment of geodesic complex curve
tangent to $(p_0;\xi_0)$ iff $F$ is K\"ahler at~$(p_0;\xi_0)$, that is iff
$$T_{\alpha i\smb\mu}(p_0;\xi_0)\,\xi_0^i\bar{\xi_0^\mu}=0,\qquad\alpha=1,
        \ldots,n.\neweq\eqdXkahler$$
The segment, if exists, is unique. Furthermore, if $F$ is complete then the
segment of geodesic complex curve actually extends to a whole geodesic complex
curve.
 
\pf One direction is known. Conversely, assume \eqdXkahler\ holds. By
Theorem~\dC,
$$G_{\alpha\sigma}\Sigma^\sigma_{;i\smb\j}\,v^i\bar{v^j}=0;$$
hence, by Proposition~\dXkahler, $X\bigl(T_{\alpha i\smb\mu}\,v^i\bar{v^\mu}
\bigr)=0$. So $T_{\alpha i\smb\mu}\,v^i\bar{v^\mu}$ is constant (and thus zero)
along the solution of \eqdsei\ tangent to~$(p_0;\xi_0)$; the assertion then
follows from Proposition~\dT.\qedn
 
As already discussed in the introduction, one of the motivations behind this
work was to find a differential description of the properties of the Kobayashi
metric in strongly convex domains. We conclude then with the following:
 
\newthm Corollary \dnove: Let $F$ be a strongly pseudoconvex smooth complete
complex Finsler metric on a manifold $M$. Assume that $\eqduquattro$ holds or,
equivalently, that $K_F\equiv-4$ and $\eqdcro$ holds. Then $F$ is the
Kobayashi metric of~$M$.
 
\pf By Theorem~\dC, in both cases the holomorphic curvature of $F$ is~$-4$.
Furthermore, by Theorem~\dotto\ $F$ is realizable. The assertion then
follows from Proposition~\usette.\qedn
 
\setref{ZZZ}
\beginsection References
 
\book A M. Abate: Iteration theory of holomorphic maps on taut manifolds!
Me\-di\-ter\-ra\-nean Press, Cosenza, 1989
 
\art AP M. Abate, G. Patrizio: Uniqueness of complex geodesics and
characterization of circular domains! Man. Math.! 74 1992 277-297
 
\book B A. Bejancu: Finsler geometry and applications! Ellis Horwood Limited,
Chichester, 1990
 
\art Bu J. Burbea: On the hessian of the Carath\'eodory metric! Rocky Mtn.
Math. J.! 8 1978 555-559
 
\art F J.J. Faran: Hermitian Finsler metrics and the Kobayashi metric! J.
Diff. Geom.! 31 1990 601-625
 
\art H M. Heins: On a class of conformal mappings! Nagoya Math. J.! 21 1962
1-60
 
\book K1 S. Kobayashi: Hyperbolic manifolds and holomorphic mappings! Dekker,
New York, 1970
 
\art K2 S. Kobayashi: Negative vector bundles and complex Finsler structures!
Nagoya Math. J.! 57 1975 153-166
 
\art L L. Lempert: La m\'etrique de Kobayashi et la repr\'esentation des
domaines sur la boule! Bull. Soc. Math. France! 109 1981 427-474
 
\book P M.Y. Pang: {\rm Finsler metrics with the properties of the Kobayashi
metric on convex domains}! Preprint, 1990
 
\book R H.L. Royden: {\rm Complex Finsler metrics}! In {\bf Contemporary
Mathematics. Proceedings of Summer Research Conference,} American Mathematical
Society, Providence, 1984, pp.~119--124
 
\book RW H.L. Royden, P.M. Wong: {\rm Carath\'eodory and Kobayashi metrics
on convex domains}! Preprint (1983)
 
\book Ru1 H. Rund: The differential geometry of Finsler spaces! Springer,
Berlin, 1959
 
\art Ru2 H. Rund: Generalized metrics on complex manifolds! Math. Nach.! 34
1967 55-77
 
\art S M. Suzuki: The intrinsic metrics on the domains in $\C^n$! Math. Rep.
Toyama Univ.! 6 1983 143-177
 
\art V E. Vesentini: Complex geodesics! Comp. Math.! 44 1981 375-394
 
\art W B. Wong: On the holomorphic sectional curvature of some intrinsic
metrics! Proc. Am. Math. Soc.! 65 1977 57-61
 
\art Wu H. Wu: A remark on holomorphic sectional curvature! Indiana Math. J.!
22 1973 1103-1108
 
\bigskip\bigskip
\line{\vbox{\hsize=7truecm
\leftline{Marco Abate}
\leftline{Dipartimento di Matematica}
\leftline{Seconda Universit\`a di Roma}
\leftline{00133 Roma, Italy}}\hfill\vbox{\hsize=7truecm
\leftline{Giorgio Patrizio}
\leftline{Dipartimento di Matematica}
\leftline{Seconda Universit\`a di Roma}
\leftline{00133 Roma, Italy}}}
\medskip
\leftline{June 1992}
\vfill
\bye